\documentstyle[12pt,psfig] {article}
\newtheorem{theo}{Theorem}
\newtheorem{lem}{Lemma}

\newtheorem{defi}{Definition}

\newtheorem{prop}{Proposition}

\newtheorem{cor}{Corollary}
\newtheorem{rem}{Remark}

\begin {document}
\bibliographystyle{plain}
\title{ Weak Limits of Riemannian Metrics in Surfaces with integral Curvature Bound}

\author{  Xiuxiong Chen\\
   Department of Mathematics \\
   Stanford University \\
   Palo Alto, CA 94305 }
\date{ Sep. 19, 1996}
\maketitle      
\section{ Introduction}
{\bf 1.1 Introduction to  the problem.} We study the limit of a sequence 
of Riemannian metrics on a  surface under some suitable conditions.
Let $\Omega $ be any open domain and let 
${\cal{G}}(\Omega)$ be the set of all smooth Riemannian metrics on 
$\;\Omega.\;$
 Any two metrics  $g_1,\;g_2 $ in ${\cal{G}} (\Omega)$ 
are called pointwise conformally equivalent if they are related under 
multiplication by a smooth, positive function on $\Omega.\;$ This relation
 is denoted by $g_1\propto  g_2.\;$ Define the curvature energy function and area 
function for a metric $g$ in $\Omega$ as follows:
\[  K(g,\Omega)  = \int_{\Omega} K_{g} ^2 \; d\;g,\qquad A(g,\Omega) = 
\int_{\Omega} d\;g ,\]
where $K_{g}$ is the scalar curvature of $g,$ and $d\,g$ is the area (volume) element.
  For a given metric 
$g_{0}$ in $\Omega$, define a function space
 $\overline {{\cal{S}} (g_{0},C_1, C_2,\Omega)}\;$ 
to be the completion of the following set under any reasonable topology:
\[ {\cal{S}} (g_{0},C_1, C_2, \Omega) = \{ g \in {\cal{G}}(\Omega)  | 
g\propto  g_0,\; A(g,\Omega) = C_1,\;K(g,\Omega) \leq C_2\},
 \]
here $ C_1 $ and $ C_2$ are generic constants.\\

We are mainly concerned with the two following questions: (a) Given a 
sequence of metrics $\;\{g_k ,\;k \in \bf N \} $ in $ {\cal{S}} 
(g_{0},C_1, C_2, \Omega),\; $ what is the set of its cluster points? (b) When can 
one conclude that there must exist at least one limit and, if so, 
what are its geometric properties?  We have constructed an example of a
sequence of metrics with no subsequence that converges in the elementary 
sense.  Therefore, one must devise a geometrically reasonable topology in 
the function space of metrics; in particular, the area functional 
is continuous and the energy functional is lower-semi-continuous. \\

\begin{figure}
\centerline{\psfig{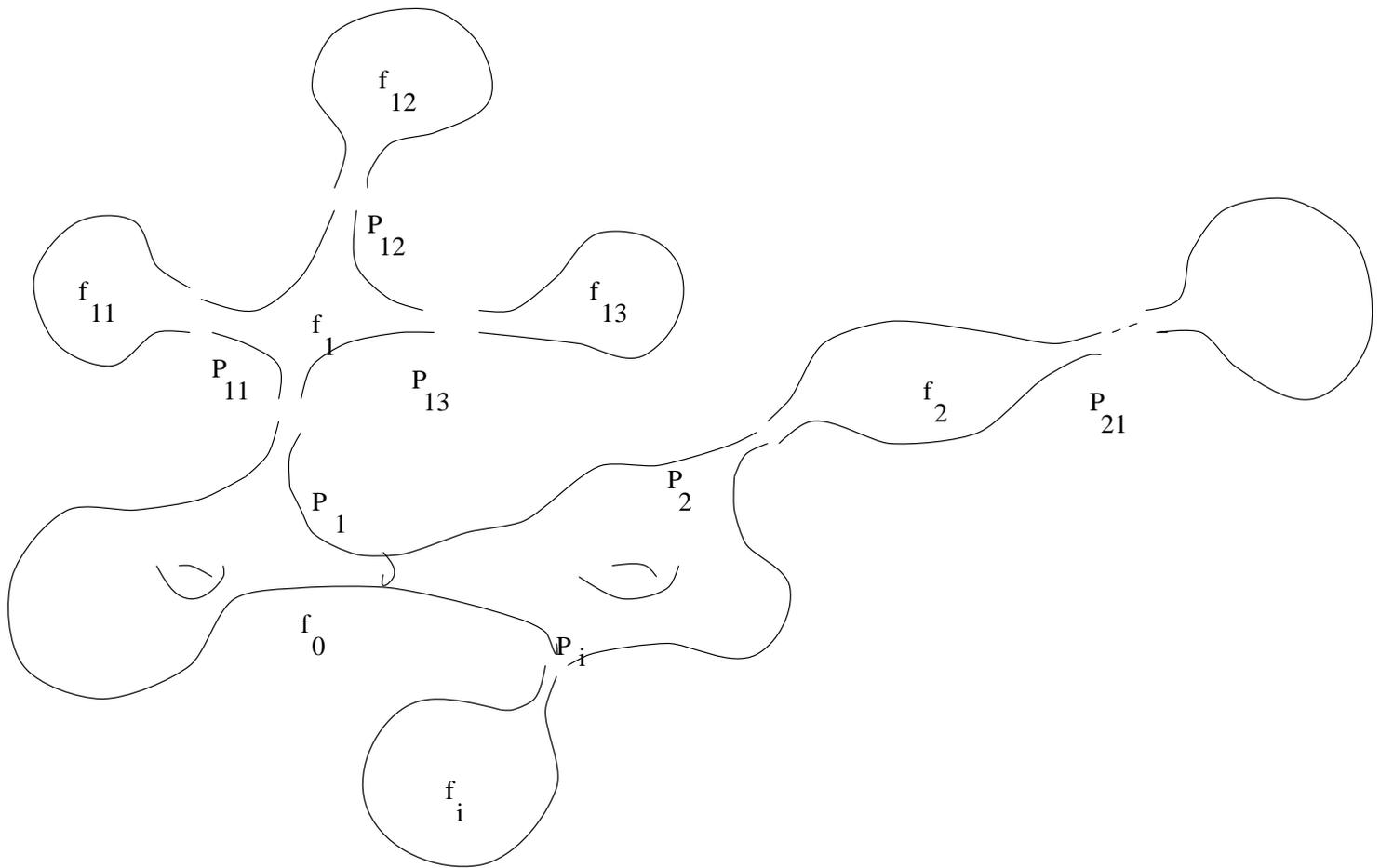}}
\caption{Bubbles on Bubbles}
\label{fg:closed}
\end{figure}

Our main result is Theorem A at section \ref{th:theorem A}.
 It could be summarized
as the following: As is shown in Figure \ref{fg:closed}, 
there is a subsequence of $\{g_n\}$
which locally weakly converges in $H^{2,2}$ (functions up to second derivative
are in $L^2$) to a Riemannian metric $f_0.\;$
However, this weak convergence is not on all of the surface $\Omega,$
but on $\Omega $ with a number of points $\{p_i\}$ deleted. There is a positive
 amount of energy and area concentrations at each point $p_i.\;$
At each point $p_i,$ we use a rescaling argument to construct 
a sequence of Riemannian metrics in $S^2$ with a small disk deleted
(the size of the disk approaches zero when the sequence takes a limit).
This renormalized sequence of metrics then (have a subsequence) converges to a
metric $f_i$  in $S^2$ with a finite number of points deleted. We then call
this metric  a ``bubble metric.''
Iterating this process at each
new bubble point of $f_i$, and so on.  The final ``limit'' of
the subsequence (passing to the diagonal subsequence) is a disjoint
 union of these bubble metrics, which are defined
in different surfaces. 
Each metric in the ``limit''  has a special
property that if it vanishes at one point in its domain, 
it then vanishes everywhere in its domain. 
 While a bubble metric might be a metric in 2-sphere with constant  curvature,
 generically it is a metric defined on a punctured sphere  and  it 
has a singular angle at each punctured point.\\

\noindent {\bf 1.2 Extremal K\"ahler metrics.} The proposed problem
is  motivated from the study of the existence problem of extremal K\"ahler
metric in a K\"ahler manifold $M.\;$ An extremal K\"ahler 
metric is a critical point for the energy functional: 
\[  E(g) = \int_M \;K_g^2 \;d\,g \]
on the space of K\"ahler metrics in a fixed K\"ahler class. The Euler
equation for the critical metric is (assume $\partial M =\emptyset$):
\[   {K_g}_{,\alpha \beta} = 0,\qquad \forall\;  1\leq \alpha,\;\beta \leq n,                    \]
where $n$ is the complex dimension of the manifold. \\
 
This problem first appeared in \cite{calabi82}.  
 E. Calabi proposed to use the heat flow
method to solve this problem when the manifold 
admits no holomorphic vector field. The heat flow he suggested is: 
\[   {{\partial g_{\alpha \overline{\beta}}}\over{\partial t}} = 
 {K_g}_{,\alpha \overline{\beta}},
\qquad\forall\;  1\leq \alpha,\;\beta \leq n.  \]
 This flow indeed decreases the energy function $E$ along its trajectory.
 To obtain the long term existence and the convergence as $ t\rightarrow \infty,\;$
 one needs to understand the following question: 
what is the weak topology of the  set of metrics in a 
fixed K\"ahler class with bounded energy? 
 An initial approach to this question would be 
to study it in the case of complex dimension one, reducing the problem to 
the one just described. Observe that any Riemannian metric on a surface is
also a K\"ahler metric; any two metrics $g_1,\; g_2$ are in the same K\"ahler
class if and only if $ g_1 \propto g_2\;$ and $ \;\int_\Omega d\;g_1 =\int_\Omega d\;g_2$ if $ \partial \Omega = \emptyset. $\\ 

\noindent {\bf 1.3 Uniformization theorem and Dirichlet Problem.} 
The selection of $L^2$ norm (rather than any $L^p$ norm with $p > 1$)
 of the  scalar curvature as energy function is not essential 
 as far as the weak topology of the function space is concerned.
 It is significant, however, if we 
consider the corresponding variational problem. The  Euler equation
  of  the energy  functional  is: 
\begin{equation} \triangle_g \;K_g + K_g^2 = C \; {\rm (generic\; constant).}
\label{eq:real:euler}
\end{equation}

This Equation is called the extremal equation. Any metric satisfies
equation (~\ref{eq:real:euler}) is called an extremal metric, even if it
is only a stationary point of the energy functional. \\

Let $ \Omega $ be any domain with smooth boundary; let $g_0$ be a smooth 
metric in $\Omega$ which could be extended smoothly to a slight larger domain.
We want to ask if there always exists a metric $g,$ in a pointwise 
conformal class of $g_0,$ which solves 
equation (~\ref{eq:real:euler}) and satisfies the Dirichlet boundary condition:
\begin{equation}
  g|_{\partial \Omega} = g_0|_{\partial \Omega},
\qquad {\partial g \over \partial n}|_{\partial \Omega}
 ={\partial g_0 \over \partial n}|_{\partial \Omega}. \label{eq:direchlt} 
\end{equation}
\noindent {\bf Conjecture 1. } {\it There always exists a solution to equation (1) with
Dirichlet boundary condition (\ref{eq:direchlt}), while solution metric
is pointwise conformal to the initial metric $g_0$.}\\

 The Euler equation (~\ref{eq:real:euler}) has an equivalent complex version:
\begin{equation}
 {\partial \over \partial {\overline{z}}} {K_g}_{,zz} = 0,\label{eq:complex:euler}
\qquad  {K_g}_{,zz} = {{\partial^2 K_g }\over{\partial z^2}} - 2\cdot
 {{\partial K_g }\over{\partial z }} \cdot {{\partial \varphi}\over{\partial z}},
\end{equation}
 where  $g = e^{2\varphi} |d\,z|^2 $ locally.\\

 The Euler equation has two important special cases:
 the first special case is the following
\begin{equation}
{K_g}_{,zz} = 0,
\label{eq:hciu}
\end{equation}
while the second special case is the following
\begin{equation} {K_g} \equiv C,\qquad {\rm or}\; - \triangle \varphi = C \cdot e^{2 \varphi}.
  \label{eq:constant}
\end{equation}

Any metric solves the equation (~\ref{eq:hciu})
 has a special property that the Hessian of its
curvature is proportional to the metric tensor.  Therefore, we may denote these
metrics as HCMU metrics (``Hessian of Curvature of Metric is Umbilical''). 
 If the Conjecture 1
were established, it would be desirable to understand the obstructions for
the existence of any HCMU metric and  obstructions to the existence of
any constant curvature metric in a domain with appropriate Dirichlet boundary
condition~(\ref{eq:direchlt}). \\

In the special case when $\partial \Omega = \emptyset$, any extremal metric also
has a constant curvature.  Recall that the classical uniformization theorem 
in a surface with no boundary asserts that any Riemannian metric is pointwise
conformal to a metric with constant curvature.  Therefore, the Conjecture 1 
(if proved), would generalize  the classical uniformization
theorem in a surface with no boundary to any domain with smooth boundaries. \\

Consider another special case where the boundaries are a set of isolated points.
To replace the Dirichlet boundary conditions, one  requires all of the metrics
 have a prescribed conical angle at each boundary point. 
Such a surface is called ``a surface with conical singularities''
 (see \cite{MTrojanov:conical91} for definition).\\

\noindent {\bf Open Problem 1.} {\it Is any  Riemannian metric on a singular surface 
 pointwise conformal to an extremal metric with the same 
angle at each singular point.}\\

In this special case, there have been plenty of attempts
 ( mostly by analysts) 
to generalize the classical uniformization theorem to  surfaces with
conical singularities. Most work has  concentrated on finding  a
constant curvature metric in a pointwise conformal class.
 However, we believe our approach may be more fruitful,
 since not all surfaces with conical singularities support
a  constant scalar curvature metric.  Our program involves two related
but independent problems. The first problem is to use direct variational
 method to give a positive answer to the above problem. For this purpose,
 we need to study the weak compactness of the function space of Riemannian metrics with
finite energy and area (which is the subject of this study).
 The second problem is to study the obstructions of 
existence of any HCMU metric 
and constant curvature metric in such surfaces.  This second problem
is discussed in \cite{chen943}, where we give a necessary condition 
for these surfaces to admit any HCMU
metric with non-constant curvature.\\

\noindent {\bf 1.4 Bubbling phenomenon.} An important feature
of  Theorem A is ``bubbling phenomenon.'' 
The bubbling phenomenon was first observed by
 Sacks-Uhlenbeck \cite{Uhlenbeck} in 1979,
 when they studied the existence theorem for harmonic maps
between two spheres. Since then, it has been studied and recognized in a wide
variety of geometric differential equations (see \cite{ParkerWolfson} for
further references). The solution spaces to these
equations are non-compact in any reasonable topology. The key
observation was that the non-compactness is
associated with the concentration of the energy density at isolated
points and that, by using the conformal invariance
of the equations, one could renormalized the solutions around
these points to obtain other solutions. This re-normalization process
is commonly referred to as ``bubbling.''\\

Our ``bubbling'' procedure appears to be 
  very similar to the re-normalization process employed first
by  Sacks-Uhlenbeck in 1979.   However, there are some significant differences.
First, the function space is not a solution space
of any elliptic equation. Second, in  most geometric
problems where the bubbling phenomenon occurs, 
the energy function involves only the first derivatives
of the ``function'' in the solution space.  However, the energy
 functional here involves the second derivatives. To make the matter worse,
it involves only the Laplacian of the conformal parameter function, which
exerts a  very weak control on the size of the metric.
These differences dictate  a new approach
other than the standard one to solve the problem.
For instance, in most of these problems where bubbling
phenomenon occurs, one usually obtains a weak convergence result without too much 
difficulty. The hard part is to show that the bubble points are
isolated. However, we have to do it exactly in the opposite order here.
 The definition of a ``bubble point'' then becomes rather tricky, because
there is no  weak convergent subsequence to work with. To overcome
this difficulty,  we introduce
the notion of ``pseudo bubble point,'' where a subsequence of 
metrics has a positive amount of energy and area concentration.
Unfortunately, the set of pseudo bubble points could be a
 dense set in the domain.\\
 
\noindent {\bf 1.5 Thick-thin Decomposition.} 
 In a higher dimensional
compact manifold, the  Cheeger-Gromov theorem \cite{cheeger86} states
that  any sequence of metrics in a compact manifold has a convergent
subsequence,  provided that the  sectional  curvature is uniformly bounded,
the volume is bounded from below, and the diameter is bounded from above. 
Similar results to \cite{cheeger86}  were obtained 
in \cite{Gromove81},\cite{Green88} and \cite{Peter87} 
as well.
 The following corollary of the theorem A  could be regarded as a 2 dimensional version
 of the Cheeger-Gromov thick-thin decomposition theorem, 
under a weaker integral condition on
the curvature tensors.
\\

\noindent {\bf Corollary B.} {\it  For any locally weakly convergent sequence
 of surfaces $ \{(\Omega,g_n),\;n \in {\bf N}\}\;$
where  $ g_n \in {\cal{S}}(g_{0}, C_1, C_2, \Omega), $ and for  any number $\epsilon > 0,$
 there exist two integers
$ N_{thick} \;{\rm and }\; N_{thin} $ which depend only on $\epsilon$ and the total
energy $\sqrt{C_1 \cdot C_2} $ of this sequence (independent of $n$).
There exists a decomposition of
$(\Omega,g_n)$ into $ N_{thick}$ of thick components 
$ \{(\Omega_{\alpha},g_n|_{\Omega_{\alpha}})\} $ (indexed by $ I_{thick}$) and $N_{thin}$ of 
thin components $\{(\Omega_{\beta},g_n|_{\Omega_{\beta}})\}$ (indexed by $ I_{thin}$), 
such that (see Figure~\ref{fg: thin-thick} on p.~\pageref{fg: thin-thick}): 1) $ \Omega = \displaystyle {\bigcup_{\alpha \in I_{thick}}} \Omega_{\alpha} 
\displaystyle {\bigcup_{ \beta \in I_{thin}}} \Omega_{\beta}. \;$ 2)
For any fixed $\alpha \in I_{thick},$ except one, $ (\Omega_{\alpha}, g_n|_{\Omega_{\alpha}}) $
locally weakly converges to a metric in $S^2$ with a finite number of 
small disks deleted; the other thick component locally weakly converges to a metric
 in $\Omega$ with a finite number of disks deleted.  Moreover, the size of each
 deleted disk could be made as small as needed.  3) Each thick component is 
self-connected; however, they are mutually disconnected if all of
the thin components are removed from the surface.
4) Each thin component is topologically  $ S^1 \times (a,b)$ and  the
length of any concentric circle $S^1 \times \{x\} ( x \in (a,b)) $ is  strictly
less than $\epsilon.$ }\\
 
\noindent {\bf Remark 1.} {\it The terms ``thin'' and ``thick''  used here,
 are slightly
different from  what are originally used in \cite{cheeger86}.
 For instance, the metrics in
a  thin part in above corollary do not necessary have a 
lower bound on the scalar curvature. }\\
 
\begin{figure}
\centerline{\psfig{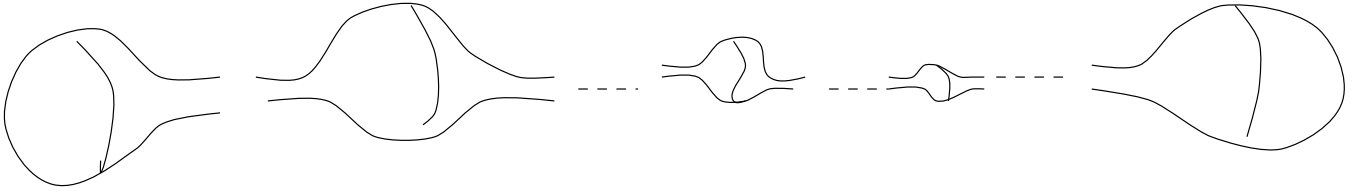}}
\caption{Rotationally symmetric Bubbles}
\label{fg:mid-bubbles}
\end{figure}

We initially hoped that both numbers $N_{thin}$ and $N_{thick}$ would be independent
 of $\epsilon.\;$ However, we have constructed a sequence of rotationally
symmetric metrics in ${\cal{S}}(g_0,C_1, C_2, S^2)$ such that this sequence yields
as many thick components as needed when $\epsilon \rightarrow 0,$ without
incurring a blowing up of the energy functional. We first constructs a sequence of
metrics in a  sequence of disks where the boundary of each disk is a smooth closed geodesic;
the length of the boundary geodesic tends to 0, while both the energy functional
and area are kept uniformly bounded from above 
(see example 2 in p. \pageref{example} for details). We then construct a sequence of
metrics in a sequence of cylinder where both boundary circles are geodesics;
the length of the boundary geodesics tends to zero, while the energy functional
and area functional could make to be arbitrarily small. Using these metrics as building
block, we could construct a sequence of metrics with
bounded energy and area as in Figure \ref{fg:mid-bubbles}, where the limit
of metrics splits into as  many parts as desired.  Henceforth, Corollary {\bf B} 
in its present form is the best one we could expect.\\

Motivated by the work of \cite{cheeger86}, C. Barvard and P. Pansu \cite{BavardP88}
 studied the divergence problem of a sequence of metrics in any surface with
 pointwise curvature  bounded, allowing the conformal structure to be varied. 
They have constructed some examples which  show that the compactness  fails 
if the conformal structure is not bounded. As a matter of fact, the weak 
compactness still fails even
 if the conformal structure is fixed.  Following the work of
C. Barvard and P. Pansu \cite{BavardP88}, M. Trojanov 
\cite{MTrojanov:concentration91} first considered a sequence 
of Riemannian metrics in a surface with a $L^p (\;\forall \;p > 1)$ 
norm of curvature (with respect to a fixed background metric) uniformly
 bounded from above. He then showed such a sequence of metrics 
either has a convergent subsequence or has at least one singular point.
However, the ``bubbles on bubbles'' phenomenon is not observed
 in \cite{MTrojanov:concentration91}\\

\noindent {\bf 1.6 Analytical approach.} In a local coordinate
system,  any Riemannian metric  $g$ could be expressed in
terms of its conformal parameter function $\varphi:$
\[ g = e^{2\varphi} (d\,x^2 + d\,y^2).\]
Therefore,  $g$  can then be regarded as a solution 
of  scalar curvature equation:
\begin{equation} -\triangle\;\varphi = K  \cdot e^{2\varphi}.
 \label{eq:anaEQ} \end{equation}
 
H. Brezis and F. Merle \cite{BrezisM91} 
 had studied the weak compactness
of the solution space of this equation.  They consider a sequence 
of pointwise conformal metrics in an  open disk. It is assumed that
the $L^p (\forall \;p>1) $  norm of curvature is
 uniformly  bounded from above  and the curvature function
is  non-negative.  They  \cite{BrezisM91} 
observed only the  first level of bubbles, but not bubbles 
 on bubbles phenomenon. \\

In both problems, difficulties arise because the right side of 
equation (\ref{eq:anaEQ}) is only in $L^1.\;$ Interested readers are encouraged
to compare the main theorems of \cite{BrezisM91} with the 
 Theorem~\ref{theo:lo:weak} and~\ref{theo:lo:bubble}
in Section 3, where  the problem is discussed from an analytic perspective.
 There are some striking similarities which underscore the  connections
between theses two problems. However, there are also some fundamental differences
 between these two problems. It is assumed in \cite{BrezisM91}
  that either the scalar curvature function 
is non-negative, or the area element is in $L^{p'}
 ({1 \over p} + { 1 \over p'} = 1).\;$ The compactness fails in our problem  precisely
because that the scalar curvature function changes sign and the area element
is only in $L^1.\;$   \\

\noindent {\bf 1.7 Organization}. In Section 2, we introduce the corresponding
 local weak compactness
problem and conclude a local version of weak convergence theorem. Also in this
section, we analyze the sequence of metrics near a bubble point via blowing
up and conclude a theorem of bubbles on bubbles. This Section is the central 
piece of this work. In Section 3, we essentially restate the weak compactness
theorem in a geometric context. In Section 4, we outline a bubbling procedure
and obtain a theorem of bubbles on bubbles.

\section { Local problem from an analytic viewpoint}

\subsection{Introduction}
In this section, we consider the problem of weak compactness of a sequence of metrics
 in a local coordinate disk. In one coordinate chart $(D,z)$, any
 metric $g$ is defined as:
\begin{equation} g = e^{2 \varphi} (d x^2 + d y^2),  \label{eq:lo:metric}
\end{equation}
and  the  curvature function is:
\begin{equation} K = - { {\triangle \varphi}\over {e^{2 \varphi}}}. 
\label{eq:lo:curv}
 \end{equation}

A metric $g$ is said to have a finite area $C_1$ and a finite energy $C_2$ 
 if and only if the following conditions are met:
\begin{equation}
 \left\{ \begin{array}{ccc} \int\limits_{D}
 e^{2 \varphi} dx dy & \leq & C_1, \\ 

     \int\limits_{D} {{(\triangle \varphi) ^2} 
\over {e^{ 2 \varphi}}} d\, x d\,y & \leq  & C_2.
\end{array} \right. 
\label{eq:lo:bound}
\end{equation}

A sequence of metrics $\{g_n\}$ where $g_n = e^{ 2\varphi_n} (d\,x^2 + d\,y^2)$
is said to have finite area $C_1$ and energy $C_2$ if and only if each $\varphi_n$
satisfies the inequality~(\ref{eq:lo:bound}). From this point on, in this
Section, we will use either $\{\varphi_n\}$ or $\{g_n\}$ to denote
a sequence of metrics with finite area $C_1$ and energy $C_2,$ unless
otherwise specified. \\

 The questions raised in Section 1.1 are : (1)
for a sequence of metrics $\{\varphi_n \}$
 satisfies the inequalities (\ref{eq:lo:bound}),
does this sequence of functions have a  uniform bound in $L^{\infty}(D)$?
(2) what is the weak limit of $\{\varphi_n\}$ under some reasonable topology?
\\

{\small
\begin{rem} H. Brezis and F. Merle\cite{BrezisM91} considered a sequence of metrics
 $\{\varphi_n\}$ satisfies the following equation:
\[ - \triangle \varphi_n = K_n \cdot e^{2 \varphi_n}\]
in an open disk $D$, where $ K_n \geq 0 $ and $K_n \in L^{p}(D),\;e^{2 \varphi_n}
\in L^{p'}(D) $ where $ {1\over p} + {1 \over p'} = 1).\;$  They proved that
one of the following three alternatives holds true (mutually exclusive):
\begin{enumerate}
\item Vanishing case: $ \varphi_n \rightarrow - \infty $ uniformly in any compact
subset of $D.\;$
\item Convergence: there exists a function $ \varphi \in H^{2,2}(D) $ such
    that $ \varphi_n \rightharpoonup  \varphi $ weakly in $H^{2,2}_{loc}(D). $
\item There exists a finite number of bubble points $\{p_1,p_2,\cdots,
p_m\}$ such that as a measure, 
\[ K_n \cdot e^{2\varphi_n} \rightharpoonup \sum_i \alpha_i \cdot \delta_{p_i}.  \] 
\end{enumerate}
 
They conjectured that $ \alpha_i = 4 \pi \cdot m_i$ for some  integer $m_i.\;$
This conjecture  was proved by Y. Y. Li and I. Shafrir\cite{LiSh94}. 
However, it  remains open whether
$m_i$ actually equals $1$.
 
Our problem differs from the problem consider by H. Brezis and F. Merle
 significantly. We quote their results here for comparison. The non-compactness
 occurs in our case is precisely  because the curvature
changes sign in a small neighborhood and the area element is not in
$L^{p'}$ for any $p'> 1.\;$
\end{rem}
}

For any sub-domain $\Omega$ in $D$, re-label the 
energy and area for a conformal parameter functions as:
\[
 A_c(\varphi, \Omega) = \int_{\Omega} e^{2 \varphi} d\,x d\,y, \qquad
 K_c(\varphi, \Omega) = \int_{\Omega} {{(\triangle \varphi) ^2} \over {e^{ 2 \varphi}}}
  d\,x d\,y.
\]

 A ``$0$''  metric should have ``$0$''  area and energy. Since a ``$0$''
metric has a conformal parameters function $-\infty$,  we define:
$A_c(-\infty,\Omega) = K_c(-\infty,\Omega)= 0.\;$\\

 For the convenience of notations, we add ``$-\infty$'' into
 $H^{2,2}(\Omega)$.  The resulted space is denoted by 
${\hat{H}}^{2,2}(\Omega).\;$
A sequence of functions $\{\varphi_n\} \in H^{2,2}(\Omega) $
 weak converges to a function $\varphi_0$ in
 ${\hat{H}}^{2,2}_{loc}(\Omega)$ if one of the
following two alternatives holds true (mutually exclusive):
\begin{enumerate}
\item (Vanishing case): If $\varphi_0 \equiv -\infty, $ 
then $\varphi_n \rightarrow -\infty$
 uniformly in any compact sub-domain of $\Omega.\; $ 
\item (Non-vanishing case): If $\varphi_0 \in H^{2,2}(\Omega), $ then
$ \varphi_n \rightharpoonup \varphi_0\;{\rm weakly\; in} \; H^{2,2}_{loc}(\Omega).$
\end{enumerate}  

\begin{defi} A point $p$ is said to be a bubble point of $\{\varphi_i\}$ if
for any $r > 0$,
\begin{equation}
 \displaystyle {\underline{\lim}}_{n \rightarrow \infty}
\int\limits_{D_{r}(p)} {{(\triangle \varphi_n) ^2} \over {e^{ 2 \varphi_n}}}
 dx\,dy \geq \alpha>0,\qquad 
 \displaystyle {\underline{\lim}}_{n \rightarrow \infty} \int\limits_{D_{r}(p)} 
e^{2 \varphi_n} d\,x d\,y \geq \beta>0. \label{eq:bubble}
 \end{equation}
where $ D_{r}(p)$ denotes a 
coordinate disk centered at $p$ with radius $r.\;$
The largest possible numbers $\alpha $ and $\beta$ are the concentration weights
 of the energy function
and area function at this point $p.\;$ 
\end{defi}

Clearly, if $p$ is a bubble point of $\{\varphi_n\}$, then $p$ is
a bubble point of any subsequence of $\{\varphi_n\}.$ \\

\noindent {\bf Example 1.} {\it  Let 
$g_n = {{n^2}\over{(1+ n^2\cdot |z + n^{-0.33} |^2)^2}} |d\,z|^2$ be a 
sequence of metrics in $S^2$ with a constant curvature of 1. This sequence 
of metrics then converges to $0$ at every point (including the point $z = 0$)
  on $S^2.\;$ However, the concentrations of energy and area  at $z=0$ are $4\pi, 4\pi$.
The metrics could be renormalized as: $\tilde{g}_n(z) =g_n((z-n^{-0.33})/n).\; $
This new sequence  $\tilde{g}_n$
weakly converges to a metric in $S^2$ with constant curvature.  }\\

The main theorems in this Section are:

\begin{theo} Let $\{\varphi_n,\; n\in {\bf N}\}$ be a sequence of metrics
in $H^{2,2}(D)$  with a finite area $C_1$ and energy $C_2$. There exists a 
subsequence $\{\varphi_{n_j},j \in {\bf N}\}$ of $\{\varphi_n\},$ 
a finite number of bubble points
 $\{p_1,p_2,\cdots,p_m\}(0\leq  m \leq \sqrt{\frac{C_1 \cdot C_2}{4\pi^2}})$
 with respect to $\{ \varphi_{n_j}, j \in {\bf N}\}, $ and a metric $\varphi_0
\in {\hat{H}}^{2,2}_{loc}(D \setminus \{p_1,p_2,\cdots,p_m\})$  such that:
 
\[ \varphi_{n_j} \rightharpoonup \varphi_0\;{\rm  in}\; {\hat{H}}^{2,2}_{loc}
(D \setminus \{p_1,p_2,\cdots,p_m\}).\]

 If the energy and area concentrations in each bubble point $p_i$ are
$A_{p_i}$ and $K_{p_i}$ for any $i \in [1,m],$ then:
\begin{eqnarray}
 \displaystyle{\lim_{j\rightarrow \infty}} A_c(\varphi_{n_j},D) &  
= & A_c(\varphi_0,D\setminus\{p_1,p_2,\cdots,p_m\} ) +
 \sum_{i=1}^{m} A_{p_i} \label{eq:lo:weak:area} \\
  \displaystyle{\lim_{j\rightarrow \infty}} K_c(\varphi_{n_j},D) 
& \geq & K_c(\varphi_0,D\setminus\{p_1,p_2,\cdots,p_m\} ) +
  \sum_{i=1}^{m} K_{p_i}. \label{eq:lo:weak:energy}\end{eqnarray}
\label{theo:lo:weak}
\end{theo}

\begin{rem}
The equality in formula~\ref{eq:lo:weak:energy} holds if
$\{\varphi_n\}$ minimizes the energy function. 
\end{rem}

\begin{theo}
 For any metric $\varphi$  with a finite area $C_1$ and energy $C_2$ 
 in $D\setminus \{0\}$,
  define $\phi(r) = {1 \over 2 \pi} \int_{0}^{2\pi}
 \varphi(r\cos \theta, r\sin \theta)\, d\theta.\; $
The following three statements hold true:
\begin{enumerate}
\item  $\displaystyle{\lim_{r \rightarrow 0}} (\varphi(r\cos \theta, r\sin \theta) + \ln r) = -\infty. $
\item  $\displaystyle {\lim_{r\rightarrow 0}} \;\phi'_r(r)\cdot r $
 exists and is finite.
\item There exists a constant $\beta \in (0,1) $ and two constants 
$C_3, C_4$  such that:
\[ {1\over \beta} (\phi(r) + \ln r) + C_3 \leq \varphi(r\cos \theta, r\sin\theta) 
+ \ln r \leq \beta (\phi(r) + \ln r) + C_4.\]
\end{enumerate}
\label{theo:lo:prop}
\end{theo}

\begin{theo} (Bubbles on bubbles). 
 Let $\{\varphi_n\}$ be a sequence of metrics in $D$ with finite area
$C_1$ and finite energy $C_2$. Suppose that $ p=0$ is
the only bubble point in $D$ with area concentration $A_p$
 and energy concentration $K_p.\;$ Suppose there exists a
metric $\varphi_0 \in {\hat{H}}^{2,2}(D\setminus\{p\})$ such that
$\varphi_n \rightharpoonup \varphi_0$ in
 $ {\hat{H}}^{2,2}_{loc}(D\setminus\{p\}).\;$  
A sequence of numbers $\{\epsilon_n \searrow 0\}$ can
be chosen to re-normalize the sequence of metrics as:
$ \phi_n(x,y) = \varphi_{n}(\epsilon_{n} \cdot x,\epsilon_{n} \cdot y)
 + \ln \epsilon_{n} (\forall n \in {\bf N}).\;$
There exists  a  subsequence $\{\varphi_{n_j},j \in {\bf N}\}$
 of $\{\varphi_n\},$   a 
finite number of bubble points $\{q_1,q_2,\cdots,q_m\} 
( 0 \leq m\leq \sqrt{\frac{A_p\cdot K_p}{4\pi^2}})$
 with respect to the subsequence of metrics $\{ \phi_{n_j}\},\;$
a metric $\phi_0\in {\hat{H}}^{2,2}( S^{2}\setminus 
\{\infty, q_1,q_2,\cdots,q_m\}) $ such that:
\[ \phi_{n_j} \rightharpoonup \phi_0 \;{\rm in}\; {\hat{H}}^{2,2}_{loc} 
( S^{2}\setminus \{\infty, q_1,q_2,\cdots,q_m\}).\]

If the energy and area concentrations of $\{\phi_n\}$ 
at each point $q_i$ are $K_{q_i} $ and $A_{q_i},$ then:
\begin{eqnarray}  
 A_p & \geq &  A_c(\phi_0,S^2 \setminus\{q_1,q_2,\cdots, q_m\}) 
+ \sum_{i=1}^{m} A_{q_i}\label{eq:blow:area}\\
    K_p & \geq &  K_c(\phi_0,S^2 \setminus\{q_1,q_2,\cdots, q_m\})
 + \sum_{i=1}^{m} K_{q_i}. \label{eq:blow:energy}
\end{eqnarray}
\label{theo:lo:bubble}

 If $\phi_0 \equiv -\infty$( vanishing case),  then $m \geq 2 $
and $p\, ( z = 0) $ is a bubble point of $\{\phi_{n_j}, j \in {\bf N}\}. \;$

\end{theo}

\begin{rem}
   The difference of the left side and right side of the inequality~\ref{eq:blow:area}
represents the amount of area lost during the blowing up procedure. If this
amount is zero,  there is no area trapped in the neck.
\end{rem}

 In Subsection 3.2, we prove three important lemmas (lemma 2,4 and 6), which provide a
technical foundation for the main theorems. The proof are rather technical,
 readers are then encouraged to skip Subsection 3.2 
and read the other Subsections first. In Subsection 3.3, 
we prove a weak convergence
theorem. In Subsection 3.4,  we briefly describe the properties of the
limit metrics. In Subsection 3.5, we show that a renormalized sequence of
metrics at each bubble point will have a weak convergent subsequence.

\subsection{Small energy lemmas}
In this subsection, the notion of a pseudo bubble point is introduced.  It is
subsequently used  to prove three key lemmas: lemma 2, 4 and 6.
Lemma 2  shows that the concentration of total
energy (product of curvature energy and area) at each bubble point
must be greater than $ 4\pi^2.\;$ Thus, there are at most a finite number 
of bubble points for any subsequence of metrics.  Lemma 4  shows
that if a point is not a pseudo bubble point,
 then the sequence of metrics in a neighborhood
of that point is uniformly bounded from above. Lemma 6 
shows that in any domain, if the metrics are uniformly bounded from above, then
either the sequence of metrics approaches $0$ everywhere in its domain, or a
subsequence of these metrics weakly converges in $H^{2,2}$ in any compact sub-domain.\\

 For any $p \in D,$
a small disk center at $p$ with radius $r$ will be denoted by $D_{r}(p)$.
\[ D_{r}(p) = \{(x,y) \in D | (x-x_p)^2 + (y-y_p)^2 < r^2 \}.\]

Define local  energy and area functions with respect to any point $p \in D $ as the following:

\[\begin{array}{ccc}
 K_p(r) & = & \displaystyle{{\overline{\lim}}_{ k \rightarrow \infty}}
 { \int_{D_{r}(p)} {{(\triangle \varphi_k)}^2
         \over{e^{2 \varphi_k}}}d\,x\;d\,y },\qquad \forall r > 0, \\
A_p(r) & = & \displaystyle{{\overline{\lim} }_{ k \rightarrow \infty}}
      { \int_{D_{r}(p)} e^{2\varphi_k} d\,x\;d\,y,}\qquad \forall r > 0.  
\end{array}
\]

In this definition, the limit taken is only an upper limit, since
  it is not known whether 
 $\{\varphi_n \}$ has any weak convergent subsequence.
\begin{defi} The energy  and area concentration functions
 of a sequence of metrics $\{g_k \}$ at any point 
$\;p \in M,$  are  defined  as follows:
\[
 K_p  = \displaystyle{\lim_{r \rightarrow 0}}\, K_p(r),\qquad
 A_p  = \displaystyle{\lim_{r \rightarrow 0}}\, A_p(r).
 \]
\label{defi:energy:area}
\end{defi}

Any point $p\in D$ is called a  pseudo bubble point
 if and only if $ A_p >0 $ and $ E_p > 0.\;$ Later, we could show that
$A_p > 0$ actually implies $E_p > 0.\;$ At a pseudo bubble point, there
exists a  subsequence of $\{ \varphi_n\}$ such that this
subsequence has a positive amount of area and energy concentrations
there. If we pass to
this subsequence, the pseudo bubble point becomes a real bubble point.

\begin{prop} Let $p$ be a pseudo bubble point of a sequence of metrics
$\{\varphi_n, n \in {\bf N}\}, $ there then
 exists a subsequence of $\{ \varphi_n\}$
such that $p$ is a real bubble point with respect to this subsequence.
\label{prop:psedo-bubble}
\end{prop}
{\bf Proof.} The proof is straightforward.

\begin{defi} The waist concentration function, $ l_p(\rho,\rho_{0}), $ 
  for any 
 $0 < \rho < \rho_{0}$ is defined as:

\[
 l_p(\rho,\rho_{0}) =\displaystyle {{ \underline {\lim}}_{ n\rightarrow \infty}}\;
 \displaystyle {\min_{\rho \leq r \leq \rho_{0}} }|{\partial D_r}|_{g_{n}} =\displaystyle
 {{ \underline {\lim}}_{ n\rightarrow \infty}}\;
 \displaystyle {\min_{\rho \leq r \leq \rho_{0}} } \int_0^{2 \pi} 
e^{  \varphi_n(r \cos \theta, r \sin \theta)}\, r\, d\,\theta. \]

 \label{defi:waist}
 \end{defi}

 \begin{lem} Let $\{\varphi_n\}$ be a sequence of metrics with finite area $C_1$
and finite energy $C_2.\;$ For any $\rho_{0} > 0$, we have 
$ \displaystyle {{\lim}_{\rho\rightarrow 0}}\; l_p(\rho,\rho_{0}) = 0$.
 \label{lem:waist}
 \end{lem}

 \noindent{\bf Proof.}  If the lemma is false, then
 there  exists a number $\epsilon > 0$ 
such that: 
$\displaystyle {\lim_{\rho\rightarrow 0}} \;l_p(\rho,\rho_{0}) =  2\; \epsilon > 0.\;$  Choose $\rho$ small enough so that:
\begin{equation}
  {{\epsilon^{2}}\over{ 2 \pi}}\;\ln {{\rho_{0}}\over{\rho}} > 2 \cdot C_1.
 \end{equation} 

 Since $l_p(\rho,\rho_{0}) $ is a monotonely increasing function on its variable
 $\rho > 0,\; $ 
 \[ l_p(\rho,\rho_{0}) \geq  2 \epsilon, \; \forall\; \;0 < \rho <\rho_{0}.\]

In other words,
 \[{ \displaystyle{\underline{\lim}}_{n\rightarrow \infty}} 
 \displaystyle {\min_{\rho \leq r \leq \rho_{0}}}|{\partial D_r}|_{g_{n}}
  \geq  2\; \epsilon. \]

 Fixing the pair of numbers $\,\rho,\,\rho_{0},$ there then exists 
 a number $n_0$  such that $ |{\partial D_r}|_{g_{n_0 }} > \epsilon, \; \forall\;
 \rho \leq r \leq \rho_{0}.\;$ In a local coordinate, 
 \[ |{\partial D_r}|_{g_{n_0}} = \int_{0}^{2\pi} e^{\varphi_{n_0}} r d\,
  \theta >\epsilon, \; \forall r \in [\rho,\rho_{0}]. \]

 However,
 \begin{eqnarray*}  2 \pi \cdot C_1  & \geq & \int_{\rho}^{\rho_{0}}
 \int_{0}^{2\pi} e^{2 \varphi_{n_0}}\cdot r \cdot d\theta  d\,r 
 \cdot \int_0^{ 2\pi}\, 1 \, d\theta  \\
  & \geq & \int_{\rho}^{\rho_{0}}  
 ( \int_{0}^{2\pi} e^{ \varphi_{n_0}} d\theta)^{2} r d\,r \\
 &\geq &  \int_{\rho}^{\rho_{0}}({{\epsilon}\over{r}})^{2} r d\,r \\
 &\geq & \epsilon^{2} \ln{{\rho_{0}}\over{\rho}} > 4 \pi \cdot C_1.
\end{eqnarray*}

 The last inequality holds true because of inequality (15).
Thus, $ 2\pi > 4\pi, $ which is impossible.
 The lemma is then proved.  QED.\\

The following theorem  is a generalization of the classical
 isoparametric inequality. It is a key theorem which we
will use it over and over again.

\begin{theo}(Readers are referred to \cite{Burago80} for
further reference). Let $g$ be a
 metric in an Euclidean disk $D$  such that
$\int_D |K_g| d g < \infty.\;$ For any disk
 $D_1 \subset \subset D,$  we have:

\[ \int_{D_1}\; |K_g| d\,g \geq 2\pi -
 {(\int_{\partial D_1 }  d s_g)^2 \over {2 \int_{D_1} d g}}
 = 2 \pi - {|\partial D_1|_g^2 \over {2 \int_{D_1} d\,g}}.\]  
\label{theo:lo:ineq}
\end{theo}

 \begin{lem} Let $\{\varphi_n\}$ be a sequence of metrics with finite area $C_1$
and finite energy $C_2.\;$
  If $p$ is a bubble point of $\{\varphi_k \},\;$
 then the following inequality holds true: 
 \[  \sqrt{  K_p \cdot A_p } \geq 2 \pi.\]
 \label{lem:bubble}
 \end{lem}

 \begin{rem} (a) The best constant in the above estimate is
  $ 4 \pi.\;$ \\ 
    (b)This lemma also proves that $A_p > 0 $ if and only
 if $K_p > 0. $
\end{rem}

 This lemma implies that there are only a finite number of bubble points. 
  It can also be regarded as a ``small energy lemma.'' In other words, if 
 the total energy $ \sqrt{ K(\Omega) \cdot A(\Omega)}\; $ is small enough
  ($\leq 2 \pi$), any weak convergent subsequence of metrics does not 
 have any bubble point in any compact sub-domain of $\Omega.\;$
 \\

 \noindent{\bf Proof.} Suppose $p $ is a bubble point 
 and $ A_p > 0.\; $ Let $\;\epsilon > 0 $ 
be any small positive number.  
 Recalled that $A_p = \displaystyle{\lim}_{r \rightarrow 0} A_p(r).\;$
Since $ A_p(r)$  is a monotonely increased function on its variable $r,$  then  
${\lim}_{r\rightarrow 0} A_p(r) \geq A_p > 0.\;$ 
Choose $\rho_{0}$  and for $n$ large enough:

 \[  A_p \leq A_{\rho_0}(r) = \overline{\lim_{n\rightarrow \infty}} \int_{D_{\rho_{0}}} d\, g_{n} < (1+{\epsilon \over 2}) A(p). \]
For $n$ large enough, we have
\[
 \int_{D_{\rho_{0}}} d\, g_{n} < (1+\epsilon ) A(p)
\]

  Lemma~\ref{lem:waist} then implies: 
 \[\displaystyle {\lim_{\rho \rightarrow 0}}\; l_p(\rho,\rho_{0}) = 0,\;\; 
 \forall \; \rho_{0} > 0.\]

 For any $\epsilon > 0,$ we choose a small number $\rho_{1} < \rho_{0} $
 such that  $l_p(\rho_{1},\rho_{0}) < \epsilon.\; $
  There exists a positive  number $N $
  which depends  only  on $\epsilon $ such that 
 (after passing to a subsequence): 
 \begin{equation} \displaystyle {\min_{\rho_{1} \leq r \leq \rho_{0}}} 
 |{\partial D_r}|_{g_{n}} < 2 \epsilon,\qquad \forall\; n > N. \end{equation}

 There exists a number  $\rho_{n} \in [\rho_{1},\rho_{0}]$ such that: 
 \[ |\partial D_{\rho_{n}}|_{g_{n}} < 3\; \epsilon, \qquad \forall \rho_1\leq \rho_n \leq \rho_0. \]
 
Therefore,

 \[ A_p \leq \int_{D_{\rho_{1}}} d g_{n}\leq 
  \int_{D_{\rho_{n}}} d g_{n} \leq (1+\epsilon) A(p). \]

According to Theorem 4, we have:
\begin{eqnarray*} \int_{D_{\rho_n}} |K_{g_n}| d g_n  
& \geq & 2 \pi - {{ |\partial D_{\rho_{n}}
|_{g_{n}}^{2}}\over{ 2 \int_{D_{\rho_{n}}} d\,g_{n}}}  \\ & >  & 2 \pi - 
 {{9 \epsilon ^{2}}\over{ 2 A(p)}} > 0. 
\end{eqnarray*}

The last inequality holds for any small $\epsilon > 0.\;$  Hence, we have:
\begin{eqnarray*} 
\int_{D_{\rho_{n}}} K_{g_{n}}^{2} d\,g_{n} & \geq &
 {{ ( \int_{D_{\rho_{n}}} K_{g_{n}}  
d\,g_{n})^{2}}\over {\int_{D_{\rho_{n}}} d\,g_{n}}} \\  & > & 
{{ ( 2\pi - {{9 \epsilon ^{2}}\over{ 2 A(p)}})^{2}} \over {( 1 + \epsilon ) A(p)}}.
\end{eqnarray*}

Since $\rho_{n} < \rho_{0} $ , 

\[  \int_{D_{\rho_{0}}} K_{g_{n}}^{2} d\,g_{n} > 
{{( 2\pi - {{9 \epsilon ^{2}}\over{ 2 A(p)}})^{2}} \over {( 1 + \epsilon ) A(p)}},
\;\; \forall\; \;n > N. \]

In other words,
\[ \displaystyle {\underline {\lim}_{n\rightarrow \infty}}\;
 \int_{D_{\rho_{0}}} K_{g_{n}}^{2} d\,g_{n} > 
{{( 2\pi - {{9 \epsilon ^{2}}\over{ 2 A(p)}})^{2}} \over {( 1 + \epsilon ) A(p)}}.\]

Let $  \epsilon \rightarrow 0,\; $   then let  $ \rho_{0} \rightarrow 0$ , we have:
\[ K_p = \displaystyle{ \lim_{\rho_{0} \rightarrow 0}}\;
\displaystyle{ \underline {\lim}_{n \rightarrow \infty}}\;
 \int_{D_{\rho_{0}}} K_{g_{n}}^{2} d\,g_{n} \geq {{ 4 \pi^{2}}\over { A(p)}}. \]

The lemma is then established. QED.\\

\begin{lem} If  $ \{\varphi_n\}$ has finite area $C_1$,then
$ \displaystyle {{\underline{\lim}}_{n\rightarrow \infty}}\;
 \displaystyle {\min _{ 0\leq \rho\leq r\leq \rho_{0}}}\;
 {{1}\over{2\pi}} \int_{0}^{2\pi} \varphi_{n}
(r \cos \theta, r \sin \theta)\;
 d\,\theta$ 
is bounded from above  for any interval $[\rho,\rho_{0}]$.
\label{lem:metric:cylin}
\end{lem}

\noindent {\bf Proof.} If this lemma is false, there then exists  
$ 0 < \rho < \rho_{0}$  such that:
\[ \displaystyle {{\underline{\lim}}_{n\rightarrow \infty}}\; 
   \displaystyle {\min _{\rho\leq r\leq  \rho_{0}}}\; 
 {{1}\over{2\pi}} \int_{0}^{2\pi} \varphi_{n}
(r \cos \theta, r \sin \theta)\;
 d\,\theta= \infty .\]

Passing to a subsequence if necessary, we may assume:
\[
 {{1}\over{2 \pi}} \int_{0}^{2\pi}\varphi_{n}(r\cos \theta,r \sin \theta)
 d\theta > n, 
\; \forall\; r\in [ \rho,\rho_{0}] .
\]

A Schwartz type of inequality implies:
\[ {{1}\over{2\pi}} \int_{0}^{2\pi} e^{ 2\varphi_{n}(r\cos \theta,r \sin \theta)} d\theta 
> e^{{{1}\over{2 \pi}} \int_{0}^{2\pi} 2 \varphi_{n}(r\cos \theta,r \sin \theta) d\theta }
 > e^{2 n},\; \forall\; r\in[\rho,\rho_{0}].\]

The last inequality implies:
\[ C_1 >  {{1}\over{2\pi}}\int_{\rho}^{\rho_{0}} \int_{0}^{2\pi} 
e^{ 2\varphi_{n}(r\cos \theta,r \sin \theta)}\; r d\,r  d\theta > e^{2 n} 
\int_{\rho}^{\rho_{0}} r d\,r = e^{2n} (\rho_{0}^{2}-\rho^{2})/2 \rightarrow \infty. \]
 
This is a contradiction. This lemma is proved.  QED.\\

The following lemma shows
 that the conformal parameters $\{\varphi_k\}$ 
must have a uniform upper bound, away from the set of bubble points.

\begin{lem}  Let $\{\varphi_n\}$ be a sequence of metrics with finite area $C_1$
and finite energy $C_2.\;$ If $p$ is not a pseudo bubble point of 
$\{\varphi_k,\; k \in{\bf N}\},$ i.e.,$ A(p) =0,$  there then exists a small
 neighborhood $ {\cal {O}}(p) $ of $ p  $ and a positive constant $ C $ 
 such that $ \displaystyle { \sup_{k \in {\bf N}} }
\;\displaystyle {\sup_{ q \in { \cal {O}}(p)}} \varphi_k (q) \leq C.\;$
\label{lem:regular}
\end{lem}

\noindent {\bf Proof. } Define a new function:  
\[ A_{n}(\rho) = \int_{D_{\rho}} d g_{n} = \int_{ x^2 + y^2 \leq \rho^2}
 \,e^{ 2 \varphi_n} dx dy,\qquad \forall \; n \in {\bf N},
\;\forall\; \rho \; > \;0 .\]

Choose a small coordinate disk 
$D_{r_{0}}(0)$  so that: $ 2 \cdot C_2 \cdot A_n(r_0) < \pi^2.\;$  
If this lemma is false, we could modified the sequence of metrics
 slightly so that: $\varphi_n(p) \rightarrow \infty.\;$ We want to
draw a contradiction from this assumption. \\

For any pair of numbers $ r_{1} > r_{2},\; $ consider the following:
\begin{eqnarray*}
\lefteqn{| \int_{0}^{2\pi}
 {{\partial\,\varphi_{n}(r_1\cos \theta,r_1 \sin \theta)}\over{\partial\,r}} r_1 d\,
\theta - \int_{0}^{2\pi} 
{{\partial\,\varphi_{n}(r_2\cos \theta,r_2 \sin \theta)}\over{\partial\,r}} r_2
d\,\theta | }
 \\ &  & = |\int_{r_{2}}^{r_{1}}  \int_{0}^{2\pi} {{\partial \over {\partial r}}}
  ({{\partial \,\varphi_{n} (r\cos \theta,r \sin \theta)}\over{\partial\,r}} \cdot r)
 d\,\theta \;d\,r| =  |\int_{r_{2}}^{r_{1}} 
 \int_{0}^{2\pi} ( \varphi_n''\cdot r + \varphi_n') 
 d\,\theta \;d\,r| \\ & & = |\int_{r_{2}}^{r_{1}} \int_{0}^{2\pi} \triangle \varphi_{n}(r \cos \theta,r \sin \theta) 
r  d\,\theta \;d\,r | \\ & & \leq  (\int^{r_{1}}_{r_{2}}\int_{0}^{2\pi}
 {{(\triangle  \varphi_{n} )^{2}}\over{ e ^{ 2 \varphi_{n} }}} r d\,\theta 
\;d\,r)^{ {{1}\over{2}}} \; (\int_{r_{2}}^{r_{1}} \int_0^{2\pi} e^{ 2\varphi_{n}}
 r d\theta\; d\,r )^{{{1}\over{2}}}.
\end{eqnarray*}

Since the energy of this sequence of metrics is uniformly bounded
 from above, the previous inequality implies:
\begin{eqnarray*} | \int_{0}^{2\pi} {{\partial\,\varphi_{n}(r_{1} \cos \theta,r_1 \sin \theta)}\over{\partial\,r}}\, r_1 d\,
\theta - \int_{0}^{2\pi} {{\partial\,\varphi_{n}(r_{2}\cos \theta,r_2 \sin \theta )}\over{\partial\,r}}\, r_2
d\,\theta | & \leq  \\   C_2^{1\over 2} \cdot \; (\int_{r_{2}}^{r_{1}} \int_0^{2\pi} 
e^{ 2\varphi_{n} } r d\theta\; d\,r )^{{{1}\over{2}}}. &   \end{eqnarray*}

Fixing the number $n$, observe  that
 $\displaystyle {\lim_{r\rightarrow 0} } \int_{0}^{2\pi}
 {{\partial \varphi_{n}(r\cos \theta, r \sin \theta)}\over{\partial r}} r d\theta = 0.\;$ 
Let $ r_{2} \rightarrow 0 $ and $r_1 = r,$ the result is:

\[ | \int_{0}^{2\pi} {{\partial\,\varphi_{n}(r \cos \theta,r \sin \theta)}\over{\partial\,r}} r
d\,\theta | <  C_2^{1\over 2} \cdot \; \sqrt { A_{n}(r )}. \] 

Define $\psi_{n}(r) = {{1}\over{2\pi}} \int_{0}^{2\pi} 
\varphi_{n}(r \cos \theta,r \sin \theta ) d\,\theta,\;$ then: 
\begin{eqnarray*} 
|\psi_{n}(r) - \psi_{n}(0)| & \leq & \int_{0} ^{ r} 
 {{1}\over{2\pi}}  |\int_{0}^{2\pi}{{\partial \varphi_{n}}\over{\partial \rho}} 
 \rho d\,\theta|\; \frac{ d\,\rho}{\rho} \\ & \leq & \int_{0}^{r} {{1}\over{2\pi}}\cdot C_2^{1\over 2} \cdot \sqrt{ A_{n}(\rho)} \cdot
 {{d\,\rho}\over{\rho}}.
\end{eqnarray*}

Since $\psi_n(0) = \varphi_n(p)$, therefore ($0 < \alpha < 1 $) 
\begin{equation} |\psi_{n}(r) - \varphi_n(p)| \leq  {{1}\over{2\pi}}\cdot C_2^{1\over 2} \cdot 
 \int_{0}^{r} ( {{A_{n}(\rho)}\over{\rho^{\alpha}}})^{{{1}\over{2}}}
\cdot  {{\rho^{{{\alpha}\over{2}}}}\over{\rho}} d\,\rho. \label{eq:ineq20} \end{equation}

  Following lemma 3, $\displaystyle {\lim_{n \rightarrow \infty} } \psi_{n}(r) < \infty.\; $
 It is a contradiction if the right hand side (RHS) of the previous 
inequality (\ref{eq:ineq20}) is uniformly bounded from above since $\{\varphi_n(p)\} \rightarrow \infty.\;$ However, the (RHS) of the inequality (\ref{eq:ineq20}) is bounded according
to the next lemma (choose $\alpha = 1).\;$ The lemma is then proved. QED.

\begin{lem}  Let $\{\varphi_n\}$ be a sequence of metrics with finite area $C_1$
and finite energy $C_2.\;$ Suppose  $A_p = 0.\;$ For any small number $r > 0, $ 
there exists a positive constant $C$ and a  number $N$ 
such that ($ 0< \alpha < 2$): 
\[ {{\int_{D_{\rho}} e^{2\varphi_n} d\,xd\,y} \over{\rho^{\alpha}}} =  {{A_{n}(\rho)}\over{\rho^{\alpha}}} < C,
\label{lem:no-coni} \]
if $n$ is large enough.
\end{lem}

\noindent {\bf Proof. }Choose a small coordinate disk 
$D_{r}(0)$  so that: $ 2 \cdot C_2 \cdot A_n(r) < \pi^2 ( 2 -\alpha)^2 $
(since $A_p = 0$). Let $ C $ be any number large enough such that:
\[ {{A_{n}(r) }\over{r}^{\alpha}} < C , \;\forall \; n \in
 {\bf N}.\]

It is claimed that this lemma holds true for this constant $C.\;$
 Otherwise, there exists a number  $ \rho_{n} < r$, such that 
$  A_{n}(\rho_{n}) -  C \cdot \rho_{n}^{\alpha} >  0 $. 
Consider  the function  $ F_{n} (\rho) = A_{n}(\rho) - C \cdot \rho^{\alpha}.\;$
 We have  $ F_n(r) < 0 < F_n(\rho_{n}). $  
There  then exists an interior point $r_{n} \in (\rho_{n},r_{0})$ such that:
\[ F_{n}(r_{n}) = 0, \qquad F_{n}'(r_{n}) < 0, \]
or,
\[ A_n(r_n) = C \cdot r_n^{\alpha},\qquad  
 \int_{0}^{2\pi} e^{2\varphi_{n}} r_{n} d \theta < C \cdot \alpha
\cdot  r_{n} ^{\alpha -1}. \] 
 
Using a Schwartz type inequality, we have:
 
\[
 ( \int_{0}^{2\pi} e^{ \varphi_{n}} r_{n} d \theta)^{2} <\int_{0}^{2\pi} e^{2\varphi_{n}} r_{n} d \theta \cdot \int_0^{2\pi} \,r_n\,d\theta
<   2 \pi \cdot C\cdot  \alpha \cdot \; r_{n} ^{\alpha }.
 \] 

In other words, $\; |\partial D_{r_{n}}|_{g_{n}} ^{2} < 2 \pi \alpha \cdot
  C \cdot r_{n}^{\alpha}.\;$ Therefore,
\begin{eqnarray*} \int_{D_{r_{n}}} |K_{g_{n}}| d\;g_{n} &  > & 
2 \pi - {{|\partial D_{r_{n}}|_{g_{n}} ^{2}}\over{  2 A_{r_{n}}}} \\ & > &
   2\pi -  {{ 2\pi \alpha \cdot C \cdot r_{n}^{\alpha}} \over {  2 C \cdot r_{n}^{\alpha}}} \\
& = & \pi ( 2 - \alpha) > 0.
\end{eqnarray*}

Using a Schwartz type inequality again, we have:
 
\[ \int_{D_{r_{n}}} K_{g_{n}}^{2} d\,g_n  >
{{( \pi( 2 -\alpha))^{2}}\over{ A_{n}(r_{n})}} \geq {{( \pi( 2 -\alpha))^{2}}\over{ A_{n}(r)}}. \]

Thus,
\[C_2 \geq \int_{D_{r}} K_{g_{n}}^{2} \;d\,g_n \geq  \int_{D_{r_{n}}} K_{g_{n}}^{2} 
\;d\,g >{{ \pi^{2} (2 -\alpha)^2 }\over{ A_{n}(r)}} > 2 \cdot C_2, \]
which is a contradiction. The lemmas is then proved.  QED.

\begin{lem} Let $\{\varphi_n\}$ be a sequence of metrics with finite area $C_1$
and finite energy $C_2.\;$ Suppose $ \displaystyle { \sup_{k \in {\bf N}} }\;\;
     \displaystyle {\max_{ q \in D} }  \;\varphi_k ( q)  < C_3.\;$
Let $\Omega \subset D$ be any compact sub-domain of $D.\;$
 There exists a constant
$\beta \in (0,1)$ which depends only on $C_1, C_2, C_3,$ and the domains
$\Omega, D $ such that:
\[ \displaystyle{\sup_{\Omega}}\; \varphi_k  \leq \beta \cdot
 \displaystyle{\inf_{\Omega}}\; \varphi_k + C, \qquad \forall\; k \in {\bf N}.\]
 
In particular, either $\{\varphi_{k}\} $ vanishes everywhere
 on $D$ or there exists a constant $C$ such that: 
\[
  \displaystyle {\inf_{q \in \Omega}}\,
  {\varphi_k(q)} > - C, \qquad \forall \; k \in {\bf N}. \]
\label{lem:inf+sup}
\end{lem}
\noindent {\bf Proof.}
 The conditions in this lemma are:
\[\left\{ \begin{array}{lcl} \displaystyle {\sup_{q\in D}}\;\varphi_{n}(q) 
&  < &  C_3,\\ \int_{D} {{(\triangle \varphi_{n})^{2}}\over{e^{2\varphi_{n}}}} d\,x d\,y 
&  < & C_2, \\
\int_D e^{ 2\varphi_n} dx\,dy & < & C_1\end{array}\right.
\]

 From the first two inequalities,  we imply: 
\[ \| \triangle \varphi_{n} \|_{L^{2}(D)} =
 \int_{D} (\triangle \varphi_n)^2 d\,xd\,y  < C.  \]

 Decompose the conformal 
parameter functions  $\varphi_{n}\;$ as $ \; 
\varphi_{n} =  u_{n} + v_{n},\;$  where $u_{n},v_{n} $ satisfy the following:
\[\left\{ \begin{array}{lcl} \triangle u_{n} & = & \triangle \varphi_{n},\\
u_{n}|_{\partial{D} }&  =  & 0; \end{array}\right. \] 
and 
\[
\left\{ \begin{array}{lcl} \triangle v_{n} & = & 0,\\
v_{n}|_{\partial{D} }&  =  & \varphi_{n}|_{\partial{D}}. \end{array}\right. \] 

 Clearly $ \|u_{n}\|_{H^{2,2}(D)} < C.\; $ This implies that 
 $ \displaystyle{\max_{ p\in D}}\, |u_{n}(p)| < C,\;\; \forall \; n \in 
{\bf N}. \;$ 
Since  $\varphi_{n}$ is bounded from above by the initial assumption, 
the harmonic functions $v_n = \varphi_n - u_n $
is bounded from above. For any sub-domain $\Omega \subset D,$ there
exists a constant $\beta \in (0,1)$ such that
$ \displaystyle {\sup_{\Omega}}\,(C -{v}_{n}) \leq {1\over \beta} \cdot 
\displaystyle {\inf_{\Omega}}\; ( C-{v}_{n}).\;$ Thus,
\[
 \displaystyle {\sup_{\Omega}}\,\varphi_{n} \leq \beta \cdot 
\displaystyle {\inf_{\Omega}}\; \varphi_{n} + C.\]
 QED.\\

\subsection{Locally weakly convergence}
\begin{prop} Let $\{\varphi_k,k \in {\bf N}\}$ be a sequence of metrics
in $D$ with finite area $C_1$ and energy $C_2.\;$  
There exists at most a  finite number of  bubble points
 ( bounded by $ \sqrt{\frac{C_1 \cdot C_2}{4\pi^2}}$) in $D$
 for any subsequence of metrics 
of $\{\varphi_k \}.\;$
 Moreover, there exists a subsequence of $\{\varphi_k\}$
which has a finite number of bubble points and has no additional
pseudo bubble points in $D$.
\label{prop:finite-bubble}
\end{prop}
{\bf Proof}: We first prove that there exists at most
 a finite number of bubble points
for any sequence of metrics which satisfies
 inequality~\ref{eq:lo:bound} uniformly.
Suppose that $p_{1}, p_{2},\cdots,p_{k}\;$ are all of the  bubble points. 
On one hand, we have:
\[ \sum_{i=1}^{k} A(p_{i} ) \leq  \int_{D} e^{2 \varphi_n} d\,x d\,y \leq C_1. \]

 On the other hand, lemma~\ref{lem:bubble} implies:
\[ K(p_{i}) \geq {{4\pi^{2}}\over{A(p_{i})}}, \qquad \forall \; i = 1,2,\cdots, k. \]

The total concentrated energy of this sequence of metrics at these bubble points
 must  be less than the total amount of  energy of  this sequence of metrics. Thus,
\[ \begin{array} {ccc} C_2 & \geq & 
\int_D {{(\triangle \varphi_n)^2}\over e^{2 \varphi_n}} d\,x d\,y
 \\ & \geq & \sum_{i=1}^{k} K(p_{i})
\\ & \geq  &  \sum_{i=1}^{k} {{4\pi^{2}}\over{A(p_{i})}}  \\
&\geq & {{ 4 (k \pi)^{2}}\over{\sum_{i=1}^{k} A(p_{i})}} \geq
 {{ 4 (k \pi)^{2}}\over{C_1}} .\end{array} \]

Therefore,  $ k \leq \sqrt{ {{C_1 \cdot C_2 }\over{ 4 \pi^{2}}}}.\;$  \\

 Suppose the original sequence of
metrics has $l$ distinct bubble points and has at least one additional 
pseudo point $p.\;$ Passing to an appropriate subsequence,
 (by proposition~\ref{prop:psedo-bubble}), $p$ is then a bubble point
for this subsequence. This subsequence then has $(l+1)$ distinctive bubble
points. It is claimed that a subsequence of $\{\varphi_n\}$
can be selected so that it  has only a finite number of bubble points
 and it has no additional pseudo bubble points. Otherwise, we can  
keep passing to an appropriate
subsequence to convert any additional pseudo bubble point into a new 
bubble point. Eventually, we will obtain
 a subsequence of metrics in $D$ which has more than
 $\sqrt{ {{C_1 \cdot C_2 }\over{ 4 \pi^{2}}}}$ number of bubble
points. This is a contradiction. The proposition is then proved. QED.\\

{\bf Proof of Theorem~\ref{theo:lo:weak}.}   
Passing to a subsequence of $\{\varphi_n\}$ if necessary,
 so that $\{\varphi_n\}$ has exactly $m (\geq 0)$
 number of bubble points and  has no  additional pseudo bubble points. 
Denote these bubble points by $\{p_1,p_2,\cdots,p_m\}.\;$
Choose two compact sub-domains $D^1 $ and $D^2$ so that:
\[ \{p_1,p_2,\cdots,p_m\} \subset D^1 \subset D^2 \subset \subset D.\]

Let $\epsilon > 0 $  be small enough so that
 $\{ D_{\epsilon}(p_s), 1\leq s \leq m\}$ are disjoint 
disks in $D_1.\;$  Let $D_{i,j}$ denote the following domains (see Figure~\ref{fg:domain} below):
\[ D_{i,j} = D^i \setminus (\displaystyle{\bigcup_{1\leq s \leq m}}
 D_{{1\over 2^j} \epsilon} (p_s)),\qquad i = 1,2,\; \forall \,j \in {\bf N}.   \]

\begin{figure}
\centerline{\psfig{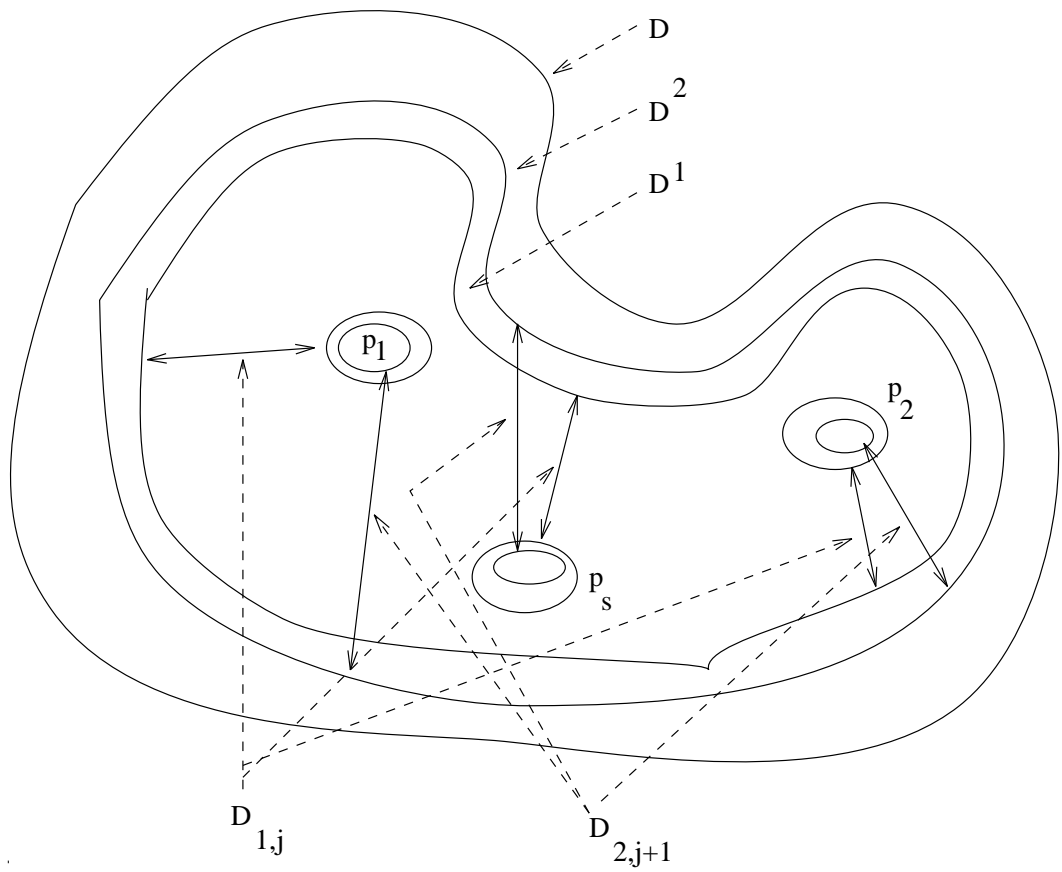}}
\caption{Compact sub-domains}
\label{fg:domain}
\end{figure}

Clearly, $D_{1,j}$ is a compact sub-domain of $D_{2,j+1}.\;$ Fixing
a number $j,$ there exists
a constant $c_j$ independent of $\{\varphi_n\}$ such that:
\begin{equation} \varphi_n(p) \leq c_j, \qquad  \forall p \in D_{2,j+1},\;\forall n \in {\bf N}.
 \label{eq:djbound} \end{equation}

If not,  there  exists a sequence of points $q_k \in D_{2,j+1} $
such that: 
\begin{equation} \displaystyle {\lim_{k\rightarrow \infty}} \varphi_k (q_k) = \infty.  
\label{eq:lo:unbound}\end{equation} 

Consider a cluster point $q \in \overline{D_{2,j+1}}$ of $\{q_k\}$
such that $ q_k \rightarrow q$ (passing
 to a subsequence of $\{q_k\}$ if necessary). 
According to the initial assumption, $q$ is not
 a pseudo bubble point of $\{ \varphi_n \}.\; $
Lemma~\ref{lem:regular} then implies that there exists a constant $C$ and
an open neighborhood $\cal{O}$ of $p$ such that
 $ \displaystyle{\sup_{n} \sup_{q\in \cal {O}}} \,\varphi_n(q) < C.\;$ 
This contradicts with equation ~(\ref{eq:lo:unbound}).
 Therefore, the inequality~(\ref{eq:djbound}) holds true. Thus, 
\[  \int_{D_{2,j}} (\triangle \varphi_n)^2 d\,xd\,y \leq e^{2 c_j} \cdot C_2.\]

According to lemma~\ref{lem:inf+sup}, 
 either $\varphi_n \rightarrow -\infty$ in $D_{2,j+1}$  or
there exists another constant $c_j'$ such that
\[ \varphi_n(p) \geq - c_j', \qquad \forall p \in D_{1,j},  \forall n \in {\bf N}.
 \]

If $\varphi_n \rightarrow -\infty$ in $D_{1,j},$ define
$\varphi_{0,j} \equiv -\infty.\;$ If $\varphi_n \not \rightarrow
-\infty $ in $D_{1,j},$ then $\varphi_n$ are uniformly bounded in
$H^{2,2}(D_{1,j}).\;$ There then exists a function $\varphi_{0,j}
\in H^{2,2}(D_{1,j})$ such that $ \varphi_n \rightharpoonup \varphi_{0,j}
$ in $ H^{2,2}_{loc}(D_{1,j}).\;$ Thus, in either case, we have:
\[ \varphi_n \rightharpoonup \varphi_{0,j} \;{\rm in}\; {\hat{H}}^{2,2}_{loc}(D_{1,j}).\]

Define $\{\varphi_{0,j}, j\in {\bf N}\}$ successively
in $D_{1,j}$ for $j=1,2,\cdots $ such that:
\[  \varphi_{n,j} \rightharpoonup \varphi_{0,j} 
\;{\rm in}\; {\hat{H}}^{2,2}_{loc}(D_{1,j}),  \]
where $\{\varphi_{n,j}\}(j > 1)$ is a subsequence of $\{\varphi_{n,j-1}\}.\;$
Consider the diagonal subsequence $\{\varphi_{n,n}\}.\;$ 
For any fixed $j>0,$ we have:
\[
\varphi_{n,n} \rightharpoonup \varphi_{0,j} \;{\rm in}\; {\hat{H}}^{2,2}_{loc}(D_{1,j}).
\]

Clearly, for any $ i > j,$ we have $ \varphi_{0,i} \equiv \varphi_{0,j}$
in $D_{1,j}.\;$ In particularly, $\varphi_{0,i} \equiv -\infty
$ if and only if $\varphi_{0,j} \equiv -\infty$ (lemma~\ref{lem:inf+sup}).
Thus, $\{\varphi_{0,j},j\in {\bf N}\}$ defines a metric
 $\varphi_0 $ in ${\hat{H}}^{2,2}(D_1 \setminus \{p_1,p_2,\cdots,p_m\})$
by
\[ \varphi_0(p) = \varphi_{0,j}(p),\qquad \forall\; p \in D_1 \setminus \{p_1,p_2,\cdots,p_m\}.            \]

Therefore,
\[ \varphi_{n,n} \rightharpoonup \varphi_0 \;{\rm in}\; {\hat{H}}^{2,2}_{loc}( D_1 \setminus \{p_1,p_2,\cdots,p_m\} ).\]                       

Denote $\{\varphi_{n,n}\}$ by $\{\varphi_n\}$ abd define
\[
A_s(j) = \displaystyle {\lim_{n \rightarrow \infty }} 
A_c(\varphi_n, D_{ {1\over{ 2^j}}\epsilon}(p_s)),\qquad K_s(j) =
  \displaystyle {\lim_{n \rightarrow \infty}} 
K(\varphi_n, D_{ {1\over{ 2^j}}\epsilon}(p_s)),\qquad \forall\, 1 \leq s \leq m.
\]

Then,
\[ A_{p_s} = \displaystyle {\lim_{j \rightarrow \infty}} A_s(j), \qquad
   K_{p_s} = \displaystyle {\lim_{j \rightarrow \infty}} K_s(j).\]

In $D_{1,j}$, we have:
\[
 \displaystyle {\lim_{n \rightarrow \infty }} A_c(\varphi_n, D_{1,j}) = A_c(\varphi_0,D_{1,j}), \qquad \displaystyle {\lim_{n \rightarrow \infty }} K_c(\varphi_n, D_{1,j}) \geq K_c(\varphi_0,D_{1,j}). 
 \]

Therefore,
\[
\begin{array}{rcc}
\displaystyle {\lim_{n \rightarrow \infty }} A_c(\varphi_n,D_1) & = & A_c(\varphi_0,D_{1,j})
+ \sum_{s=1}^{m} A_s(j), \\ \displaystyle {\lim_{n \rightarrow \infty }}
 K_c(\varphi_n,D_1) & \geq  & K_c(\varphi_0,D_{1,j})
+ \sum_{s=1}^{m} K_s(j). \end{array} 
\]

Taking limit as $ j \rightarrow \infty,$ we have
\[
\begin{array}{ccc}
\displaystyle {\lim_{n \rightarrow \infty }} A_c(\varphi_n,D_1) & = & A_c(\varphi_0,D_1\setminus \{p_1,p_2,\cdots,p_m\})
+ \sum_{s=1}^{m} A_{p_s}, \\ \displaystyle {\lim_{n \rightarrow \infty }} K_c(\varphi_n,D_1\setminus \{p_1,p_2,\cdots,p_m\}) & \geq  & K_c(\varphi_0,D_1\setminus \{p_1,p_2,\cdots,p_m\}) + \sum_{s=1}^{m} K_{p_s}. \end{array} 
\]
Let $D_1 $ and $D_2$ approach  $D$, and  use
 a similar diagonalize argument, we can show that the theorem holds
true. QED. \\

\subsection{Limit of a weak convergence sequence}

{\bf Proof of theorem~\ref{theo:lo:prop}.} We will prove 2.1, 2.2 and 2.3
separately.\\

(2.1). Let $ u = - \ln r = -\ln {\sqrt{x^2 +y^2}} $ and $ \theta = tan^{-1} {y \over x}.\;$
The domain $D\setminus \{0\}$ becomes an infinite cylinder $\{ (u, \theta) |
 0\leq u\leq \infty, -\pi \leq \theta \leq  \pi \}$ via
this transformation. Let
$\psi (u,\theta) = \varphi ( e^{-u} \cos \theta, e^{-u} \sin \theta) - u.\;$
Then $ \psi $ satisfies the following inequalities:
\begin{equation}
\left\{\begin{array}{ccc} \int_{0}^{\infty} \int_{-\pi}^{\pi}
 {(\triangle_{u,\theta} \psi)^2 \over e^{2 \psi}} d\,\theta d\,u & \leq & C_2,  
\\ \int_{0}^{\infty} \int_{-\pi}^{\pi} 
  e^{2 \psi} d\,\theta d\,u \leq & C_1, \end{array} \right.
\label{eq:lo:psi}
\end{equation}
where $ \triangle_{u,\theta} = {{\partial^2}\over{\partial u^2}} + {{\partial^2}\over{\partial \theta^2}}.\;$ To prove theorem 2.1, we only need  to show that
$\psi \rightarrow - \infty$ as $ u\rightarrow \infty.\;$ If  this
  is not true,  there then exists a positive number $C$ and a sequence
of points $\{(u_i, \theta_i), i \in {\bf N} \} (u_i \rightarrow \infty)$
such that: $\psi(u_i, \theta_i) > -C.\; $ Consider the open disk $\tilde{D} = \{ (u,\theta) | -1 < u < 1, - {\pi\over 2} < \theta < {\pi\over 2}\}.\;$ Define a new sequence of metrics in $\tilde{D}$ as:  
\[ \varphi_i(u, \theta) = \psi(u+u_i, \theta+ \theta_i), \qquad \forall i \in {\bf N},\; \forall (u,\theta) \in \tilde{D}.\]

Then $\{ \varphi_i(u,\theta), i \in{\bf N} \} $ is a sequence of functions
in $\tilde{D}$ with finite energy and area. According to theorem~\ref{theo:lo:weak}, there
exists a subsequence $\{\varphi_{n_j}, j\in {\bf N}\}$ of $\{\varphi_n\}$
and a metric $\varphi_0 \in {\hat{H}}^{2,2}_{loc}(\tilde{D}\setminus\{q_1,q_2,\cdots,q_l\})$ for some isolated singular points $\{q_1,q_2,\cdots,q_l\}$
 such that:
\[  \varphi_{n_j} \rightharpoonup \varphi_0\; {\rm in}\; {\hat{H}}^{2,2}_{loc}(\tilde{D}\setminus\{q_1,q_2,\cdots,q_l\}),\qquad l \geq 0.               \] 

The vanishing case ($\varphi_0 \equiv -\infty$) does not
occur because of
\[ \varphi_{n_j}(0,0) = \psi(u_{n_j},\theta_{n_j}) > - C, \qquad \forall j \in {\bf N}.\]

If there exists at least one bubble point $p \in \tilde{D}$, we have: 
\begin{equation} \int_{\tilde{D}} {(\triangle_{u,\theta} \varphi_{n_j})^2 \over { e^{2\varphi_{n_j}}}} d\,u d\, \theta
\cdot \int_{\tilde{D}} e^{2\varphi_{n_j}} d\,u d\, \theta > {1\over 2} E_p \cdot A_p \geq 2 \pi^2.\label{eq:polar1}
\end{equation}

If there exists no bubble point, then  $\varphi_0 \in H^{2,2}(\tilde{D})
$ and $ \varphi_{n_j} \rightharpoonup \varphi_0$
in $H^{2,2}(\tilde{D}).\;$ If $n $ is large enough, then:
\begin{equation}
 \int_{\tilde{D}} {(\triangle_{u,\theta} \varphi_{n_j})^2 \over { e^{2\varphi_{n_j}}}} d\,u d\, \theta
\cdot \int_{\tilde{D}} e^{2\varphi_{n_j}} d\,u d\, \theta > {1\over 2} \int_{\tilde{D}} {(\triangle_{u,\theta} \varphi_0)^2 \over { e^{2\varphi_0}}} d\,u d\, \theta
\cdot \int_{\tilde{D}} e^{2\varphi_0} d\,u d\, \theta > 0.
\label{eq:polar2}
\end{equation}

However, 
\[
 \int_{\tilde{D}} {(\triangle_{u,\theta} \varphi_{n_j})^2 \over { e^{2\varphi_{n_j}}}}
 d\,u d\, \theta \cdot \int_{\tilde{D}} e^{2\varphi_{n_j}} d\,u d\, \theta = 
 \int_{u_{n_j} -1}^{u_{n_j}+1} \int_{-\pi \over 2}^{\pi\over 2}
 {(\triangle_{u,\theta} \psi)^2 \over { e^{2 \psi}}} d\,u d\, \theta \cdot
\int_{u_{n_j}-1}^{u_{n_j}+1}\int_{-{\pi \over 2}}^{\pi\over 2} e^{2\psi} d\,u d\, \theta
 \]
\[ \qquad \qquad \rightarrow   0, \qquad {\rm as } \; j \rightarrow \infty. 
\]

The last formula  holds true because of inequality~(\ref{eq:lo:psi}). This
contradicts both inequalities~(\ref{eq:polar1}) and~(\ref{eq:polar2}). The
first part of the theorem is then proved. \\

(2.2). If $\displaystyle {\lim_{n \rightarrow \infty}}{ 1 \over 2 \pi}
\int_{0}^{2\pi} \varphi_r'(r \cos \theta, r \sin \theta)\, r\, d\,\theta $
does not exist, there then exist  two numbers $\alpha \neq \alpha'$
and two alternative sequence of numbers $\{\delta_i\}, \{\delta_i'\}$
such that:
\[  \delta_i < \delta_i' < \delta_{i-1} \rightarrow 0, \qquad \forall i \in {\bf N},\]
and
\[ { 1 \over 2 \pi}
\int_{0}^{2\pi} \varphi_r'(\delta_i \cos \theta, \delta_i \sin \theta) \delta_i d\,\theta
 \rightarrow \alpha,\qquad  { 1 \over 2 \pi}
\int_{0}^{2\pi} \varphi_r'(\delta_i' \cos \theta, \delta_i' \sin \theta) \delta_i' d\,\theta
 \rightarrow \alpha'. \]
Clearly,
\[ A_c(\delta_i,\delta_i') = \int_{\delta_i}^{\delta_i'} \int_{0}^{2 \pi} e^{ 2 \varphi} r d\,\theta d\,r \rightarrow 0.\]
However,
\[
\begin{array}{ccl} |\int_{\delta_i}^{\delta_i'} \int_{0}^{2 \pi} \triangle \varphi \cdot r \cdot d \, \theta d\,r | & = & |\int_{\delta_i}^{\delta_i'} \int_{0}^{2 \pi}( \varphi_r' \cdot r )'_r d\,\theta d\, r|  \\
 & = & |\int_{0}^{2\pi} \varphi_r'(\delta_i \cos \theta, \delta_i \sin \theta) \delta_i d\,\theta - \int_{0}^{2\pi} \varphi_r'(\delta_i' \cos \theta, \delta_i' \sin \theta) \delta_i' d\,\theta | \\
& \rightarrow & | \alpha - \alpha'| > 0. \end{array} 
\]

On the other hand,
\[
\begin{array}{ccc}  C_2 \geq K_c(\delta_i, \delta_i') & =& \int_{\delta_i}^{\delta_i'} \int_{0}^{2 \pi} {{(\triangle \varphi)^2}\over{ e^{ 2\varphi}}} r d\theta d\, r \\
 & \geq & {{|\int_{\delta_i}^{\delta_i'} \int_{0}^{2 \pi} \triangle \varphi \cdot r \cdot d \, \theta d\,r |^2 }\over A_c(\delta_i, \delta_i') } \rightarrow 
{{| \alpha - \alpha'|}\over {0}} = \infty. \end{array} \]

This is a contradiction. Therefore, $\displaystyle {\lim_{n \rightarrow \infty}}{ 1 \over 2 \pi}
\int_{0}^{2\pi} \varphi_r'(r \cos \theta, r \sin \theta)\, r\, d\,\theta $ does
exist. \\

(2.3). For any small $r = e^{-u} > 0, $ consider the domain $\tilde{D} = [-1,1] \times S^1.\;$ Let $\tilde{\varphi}(v,\theta) = \psi(v + u, \theta) $ (following the notations in (2.1)), then
\[
 - \triangle_{v,\theta}\; \tilde{\varphi}  =  K(v+u,\theta) \cdot e^{2 \tilde{\varphi}},\qquad \forall\, (v,\theta)\, \in \tilde{D}.
\]  

There then exists a constant $C$ such that: 
 $\tilde{\varphi} \leq C.\; $ The right hand side is bounded in
 $L^2(\tilde{D}).\;$ Let $w$ be the solution of
\[
\left\{
\begin{array}{ccl} - \triangle w & = & K(v+u,\theta) e^{2 \tilde{\varphi}},\\
                   w|_{\partial D} & = & 0.\end{array} \right. 
\]

Thus, $ ||w||_{L^{\infty}} $ is uniformly  bounded from above
(the bound is actually independent of $u$, since $L^2 $ norm of 
$ \triangle_{v,\theta}\; \tilde{\varphi} $ in $\tilde{D}$
uniformly converge to $0$ as $u \rightarrow \infty$).
The harmonic function $h = \tilde{\varphi} - w$ is bounded
from below by a constant $ - C.\;$ This follows that there exists a
constant $\beta \in (0,1)$ (independent of $u$) such that:
\[ \displaystyle{\sup_{\theta}}\,(C - h(0,\theta)) \leq {1\over \beta}
 \displaystyle{\inf_{\theta}}\,(C - h(0,\theta)),\]
Or,
\[
\displaystyle{\sup_{\theta}}\,\tilde{\varphi}(0,\theta) \leq {1\over \beta}
 \displaystyle{\inf_{\theta}}\,\tilde{\varphi}(0,\theta) + C.\]

In other words,
\[
\displaystyle{\sup_{\theta}}\,(\tilde{\varphi}(r \cos \theta, r \sin \theta)+ \ln r)
 \leq {1 \over \beta} \;
 \displaystyle{\inf_{\theta}}\,(\tilde{\varphi}(r \cos \theta, r\sin \theta)+ \ln r ) + C.
\]

Integrating on both sides over $\theta$, we obtain:
\[
{1\over \beta} (\phi(r) + \ln r) +C_3 \leq \varphi(r \cos \theta, r \sin \theta) +\ln r 
\leq \beta (\phi(r) + \ln r)+ C_4,
\]
where $C_3, C_4$ are two constants independent of $r$. QED.\\

\subsection{Bubbles on bubbles}

\begin{lem}  Suppose $D$ is a coordinate disk with radius $\rho >0 $ and 
assume \\
$g = e^{2\varphi}(d\,x^{2} + d\,y^{2})$ is a metric on $D $
with finite energy $ C_2.\;$  For any $ \epsilon > 0,$
 there then exists a constant $C_{\epsilon} > 0,$ which depends only on 
$D, \epsilon $ such that if \\
 $\displaystyle{\max_{r \leq \rho}}\;\int_{0}^{2\pi} e^{\varphi(r \cos \theta,r \sin \theta)}   r d\,\theta
 < \epsilon, $ then the following holds true: 
\[ \int_{D} e^{2 \varphi} d\,x d\,y  \leq C_{\epsilon},
\qquad {\rm and}\;\;\displaystyle \;{\lim_{\epsilon \rightarrow 0}}\; C_{\epsilon} = 0.\]
\end{lem}

\noindent {\bf Proof.} If the lemma is false, then there  exists a sequence of metrics 
 $\{\varphi_{n},\;n \in {\bf N} \} $ such that:
\begin{equation} \left\{\begin{array}{rcl} \int\limits_{D}
{{(\triangle \varphi_n)^2} \over{e^{2 \varphi_n}}} d\,x d\,y &  < & C,
 \\ \int_{D}  e^{2 \varphi_n} d\,x d\,y  & = & 1, \\
  \displaystyle {\max_{ r \leq \rho}}\;\int_{0}^{2\pi} 
e^{\varphi_n(r \cos \theta,r \sin \theta)} r d\,\theta
  &  = & \epsilon_{n} \rightarrow 0. \end{array}\right. 
\label{eq:thin}
\end{equation}

Any circle ($ |z| = \delta > 0 $) must have a zero length in the limit. 
 Therefore, $\varphi_{n}$ vanishes identically except at the origin ($z=0$). All of
the area  concentrates at the origin since the total area is fixed.  Let $\varepsilon > 0 $
 be a very small positive number and  let $ \{\delta_{n}\}$ be 
a sequence of numbers such that: 
\[ \int_{ r \leq \delta_{n} } e^{2 \varphi_{n}} d\,x\, d\,y
 = \varepsilon < {2\pi^2 \over C_2}. \] 

 Define a new sequence of metrics $ F_{n} = e^{ 2 w_{n}}(d\,x^{2} + d\,y^{2})$ as:
\[  w_{n} (z) = \varphi_{n} (\delta_{n} \cdot z) + \ln \delta_{n}, \;\forall r \leq 1.\]

For this new sequence of metrics, we have
\[\begin{array}{ccccl} \int_{D_1(0)} e^{ 2 w_{n}} & = &\int_{ r\leq \delta_{n}}\,e^{ 2 \varphi_{n}}& = & \epsilon, \\
\int_{D_1(0)} {{(\triangle w_{n})^2} \over {e^{ 2 w_{n}}}} d\,x d\,y & = &
\int_{ r\leq \delta_{n}} {{(\triangle \varphi_{n})^2} 
\over {e^{ 2 \varphi_{n}}}} d\,x d\,y & \leq & C_2.
\end{array}\]  

By theorem~\ref{theo:lo:weak}, there exists a subsequence which locally
weak converges to a metric except at a set of finite number of
 bubble points. Since each circle $|z| =\delta > 0$ has length $0$, 
 the only bubble point must be the original point $z = 0$ and all 
the area concentrated at $0$ must be less than $\varepsilon.$  
 Lemma~\ref{lem:bubble}  implies that
 total energy concentration 
 at the origin must be  bigger than ${4\pi^2 \over \varepsilon} \geq  2\cdot  C_2,$ 
which is a contradiction.  This lemma is then proved.
QED.\\

\begin{cor} Suppose  $\{\varphi_n\}$ is a sequence of metrics with finite area $C_1$
 and energy $C_2.\;$ There  exists a constant $\epsilon_{0}$  such that if 
\[\displaystyle {\max_{r\leq \rho} }\, |\partial D_r|_{g_n} =
\displaystyle {\max_{r \leq \rho}}\,\int_0^{\pi} 
e^{\varphi_{n}(r\cos\theta,r\sin\theta)} \cdot r\, d\,\theta 
\leq \epsilon_{0},\qquad \forall  n \in {\bf N}, \]

then the sequence of metrics does not have any bubble points in $D$.
\label{cor:blow}
\end{cor} 
{\bf Proof}: Let $C_{\epsilon}$ be the constant defined according to the previous
lemma. We may choose $\epsilon_0$ so small that $C_{\epsilon_0}$ satisfies:
\[  C_{\epsilon} \cdot C_2 \leq 2 \pi^2.\] 

 According to the previous lemma,  we have:
\[
 A_c(\varphi_n, D) \leq C_{\epsilon_0} \leq { 2\pi^2 \over C_2}.
\]

Thus the number of possible points $m$ must be bounded by
\[m \leq  \sqrt{ {C_{\epsilon_0} \cdot C_2 \over {4\pi^2}   }} < 1.
\]

Therefore, this sequence of metrics has no bubble points.
 The corollary is then proved.  QED.\\

\noindent {\bf Proof of theorem \ref{theo:lo:bubble}.}
Choose any small positive number $\epsilon_0 \in (0, \varepsilon).\;$  This number
 $\varepsilon\; $ serves as a scaling constant (filter).  
The sequence of functions can be modified slightly so that
 the following holds true:
\begin{equation}
  \varphi_{n}(p) =  \displaystyle {\max_{q \in D_{r_0}}} \; \varphi_{n} (q). 
\label{eq:blow0}
\end{equation}

Following  theorem~\ref{theo:lo:prop}, in a non-vanishing case, we have
$ \displaystyle {\lim_{ r \rightarrow 0 }}\; \displaystyle
{\max_{0\leq \theta \leq 2\pi}}({\varphi_0}(r \cos \theta,r \sin \theta) + \ln r )
= - \infty.\;$  There then exists a number $ r_{1} > 0 $ such that:
\[ \displaystyle {\max_{ 0\leq \theta \leq 2\pi}}  ({\varphi_0}
(r\cos \theta, r \sin \theta) + \ln r)  \ll \varepsilon, \;\;\forall\; r <  r_{1}. \]

If $ n $ is large enough, 
\begin{equation} \displaystyle{\max_{0\leq\theta\leq 2\pi}}
 ( {\varphi_{n}}(r_{1} \cos \theta,r_1 \sin \theta) + \ln r_{1})
  \ll \varepsilon, \qquad \forall \;n > N, 
\label{eq:blow1}
\end{equation}

or the length of the circle $ |z| = r_1$ is very small:
\begin{equation}
 |\partial D_{r_1} |_{g_n} = 
\int_0^{2 \pi} e^{\varphi_n(r_1\cos \theta, r_1 \sin \theta)}\, r_1\, 
d\,\theta \ll \varepsilon, \qquad n > {\bf N}.
\label{eq:blow2}
\end{equation}

According to corollary~\ref{cor:blow}, if $\varepsilon$ is small enough, 
 we can choose $\delta_{n}$ such that: 
\begin{equation} 
 |\partial D_{r} |_{g_n} = 
\int_0^{2 \pi} e^{\varphi_n(r \cos \theta, r \sin \theta)} r\, d\,\theta 
 < \varepsilon,\qquad \forall\; r_{1} \geq r \geq {\delta_{n}}, 
\label{eq:blow3} \end{equation}
and 
\begin{equation}
 |\partial D_{\delta_n} |_{g_n} = 
\int_0^{2 \pi} e^{\varphi_n(\delta_n \cos \theta, \delta_n\sin \theta)}
 \delta_n \,d\,\theta  = \varepsilon.  
\label{eq:blow4}
\end{equation}

Re-normalize this sequence of metrics as: 
\begin{equation} \phi_{n} (z) = 
\varphi_{n}(\delta_{n} \cdot z ) + \ln \delta_{n}, \qquad \forall\, |z| < {1\over \delta_n}. 
\label{eq:blow5}
\end{equation}

For any $n > 0,\; \phi_{n} $ is then defined in the disk $D_{\delta_{n}^{-1}}(0)$.
 For any fixed number $r > 0,\; \phi_{n} $ is well defined on $D_{r}(0)$ 
if $n$ is  large enough.  Moreover,  $\{\phi_{n}\}$ has a finite amount of 
energy and area since
\begin{equation}
\left\{ \begin{array} {lcl} \int\limits_{D_{\delta_{n}^{-1}}}
 {(\triangle \phi_n)^2 \over{e^{2\phi_n}} } d\,x d\,y 
  & = & \int\limits_{D_{1}}{(\triangle \varphi_n)^2 \over{e^{2\varphi_n}} }
 d\,x d\,y  \leq C_2,\\  
   \int\limits_{D_{\delta_{n}^{-1}}} e^{2 \phi_n} d\,x d\,y & = &
 \int\limits_{D_{1}} e^{ 2\varphi_n} d\, d\,y  \leq C_1.
  \end{array}\right.
\label{eq:blow6}
\end{equation}

Applying theorem~\ref{theo:lo:weak} successively to $\{\phi_n\}$ in a sequence of disks 
$D_{2^j}(0) (j = 1,2,\cdots ).\;$ In disk $D_2(0),$ there exists a subsequence
of $\{ \phi_{1n}, n \in {\bf N}\}$ of $\{\phi_n, n \in {\bf N}\},$
  a finite number of bubble points $ s_1
 = \{p_{11},p_{12},\cdots,p_{1m_1}\} (m_1 \geq 0) $ with respect to 
this subsequence, and a metric $\phi_{0,1}
\in {\hat{H}}^{2,2}_{loc} (D_2(0) \setminus \{p_{11},p_{12},\cdots,p_{1m_1}\}) $ such that:
\[
\phi_{1n} \rightharpoonup \phi_{0,1} \; {\rm in}\; {\hat{H}}^{2,2}_{loc} 
(D_2(0) \setminus \{p_{11},p_{12},\cdots,p_{1m_1}\}),
\] 

 Consider the sequence $\{\varphi_{1n}\}$
in disk $D_{2^2}(0).\;$ There exists a subsequence  $\{\phi_{2n}\}$ of $\{\phi_{1n}\}, $
  a finite number of bubble points 
$ s_2 = \{ p_{21},p_{22},\cdots, p_{2m_2}\} (m_2 \geq 0) $
 with respect to this subsequence, and a metric $\phi_{0,2} 
\in {\hat{H}}^{2,2}(D_{2^2}(0) \setminus \{p_{21},p_{22},\cdots,p_{2m_2}\}) $ such that:

\[
\phi_{2n} \rightharpoonup \phi_{0,2} \;{\rm in}\; {\hat{H}}^{2,2}_{loc} 
(D_{2^2}(0) \setminus \{p_{21},p_{22},\cdots,p_{2m_1}\}),
\] 

Clearly, the set
$ s_1 = \{p_{11},p_{12},\cdots,p_{1m_1}\} $ is a subset of 
$ s_2 = \{p_{21},p_{22},\cdots,p_{2m_1}\}$ and $ m_2 \geq m_1.\;$
Moreover, $\phi_{0,2}  = \phi_{0,1}$ when both functions are
restricted to the smaller domain $D_2$. In particular,
$ \phi_{0,1} \equiv -\infty  $ if and only if $\phi_{0,2} \equiv -\infty.\; $
\\

 In general, suppose that for any $i \leq j,$ a subsequence $\{\phi_{in}\}$ 
had been selected, and a limit metric $\{\phi_{0,i}, i\leq j\}$ had been defined
 in ${\hat{H}}^{2,2}(D_{2^i}(0)\setminus s_i) $ where $s_i$ is the set of bubble
points of $\{\phi_{in}\}$ in $D_{2^i}(0)$.  
 Consider the subsequence $\{\phi_{jn}\}$ in $D_{2^{j+1}}(0).\;$ 
There exists a subsequence $\{\varphi_{(j+1)n}\}$
of $\{\varphi_{jn}\},$  a finite number of bubble points
$s_{j+1} = \{p_{(j+1)1},p_{(j+1)2},\cdots,p_{(j+1)m_{j+1}}\}$ 
with respect to this subsequence, and a limit metric
$\phi_{0,j+1} \in {\hat{H}}^{2,2}(D_{2^{j+1}}(0) \setminus \{p_{(j+1)1},p_{(j+1)2},\cdots,p_{(j+1)m_2}\}) $  such that:
\begin{equation}
\phi_{(j+1)n} \rightharpoonup \phi_{0,j+1} \; {\rm in} \; {\hat{H}}^{2,2}_{loc} 
(D_{2^{j+1}}(0) \setminus \{p_{(j+1)1},p_{(j+1)2},\cdots,p_{(j+1)m_{j+1}}\}).
\label{eq:blow7}
\end{equation}
 
Consider the diagonal subsequence $ \{\phi_{nn}\}$. This is a subsequence 
of all the previous subsequences $\{ \phi_{jn}, n\in {\bf N}\} $ for
$j =1,2,\cdots.\;$ Therefore, all of the previous
weak convergent results hold true for this subsequence. 
In particularly, the following three  statements (for any $ j > i \geq 1$)
hold true:
\begin{enumerate}
\item  $ s_i \subset s_j.$
\item  $ \phi_{0,i} \equiv -\infty $
       if and only if $ \phi_{0,j} \equiv -\infty.$
\item  $\phi_{0,j} |_{D_{2^i}} = \phi_{0,i} $ if neither of two metrics vanishes.
\end{enumerate}

Following proposition 2, $m_j \leq \sqrt{\frac{C_1\cdot C_2}{4\pi^2}} (\forall j).\;$
 There then exists a number $N$ such that $ s_j = s_N, \forall j > N.\;$ 
We may assume that set of bubble points is: $s_N = \{q_1, q_2,\cdots,q_m\}
= \bigcup_{j} s_j (m = m_N \geq 0).\;$ \\

Define a function 
 $\phi_0 $ in $S^2 \setminus\{\infty,q_1,q_2,\cdots,q_m\}$ by:
\[ \phi_0 (p)  = \phi_{0,j}(p), \qquad \forall p \in  
D_{2^j}(0) \setminus\{q_1,q_2,\cdots,q_m\}.\]

Thus, $ \phi_0 \in {\hat{H}}^{2,2}_{loc}(S^2 \setminus \{\infty,p_1,p_2,\cdots,p_m\})$
 and the following statement holds true:
\[ \phi_{nn} \rightharpoonup \phi_0 \;{\rm in} \; {\hat{H}}^{2,2}_{loc}
(S^2 \setminus \{\infty,q_1,q_2,\cdots,q_m\}).\]

For simplicity, we re-label $\{\phi_{nn}\}$ as $\{\phi_n\}.\;$
Let $r$ be any  number large enough so that 
$\{q_1,q_2,\cdots,q_m\} \subset D_{r}(0).\;$ Consider the sequence of functions
$ \{\phi_n\} $ in $D_{r}(0)$. Suppose the concentrations of area and energy in
 the bubble point $q_i$ are $A_i$ and $ K_i.\;$ According to
 theorem~\ref{theo:lo:weak}, we have:
\begin{equation}
 \displaystyle{\lim_{n\rightarrow \infty}} \int_{D_{r}(0)} e^{ 2 \phi_n} \, d\,x \,d\,y =
  \int_{D_{r}(0)}  e^{ 2 \phi_0}\,d\,x\, d\,y + \sum_{i=1}^{m} A_i.
  \label{eq:blow:bubb1}
\end{equation}

On the other hand, 
\begin{equation}
\int_{D_{r}} e^{ 2 \phi_n}  \,d\,x\, d\,y =
\int_{D_{ \delta_n \cdot r}} e^{ 2 \varphi_n} \,d\,x\, d\,y 
 \label{eq:blow:bubb2}
\end{equation}
 
Choose a sequence of numbers $\{ \epsilon_i \searrow 0, i \in {\bf N}\}.\;$
According to the proof of proposition 1, we may have (passing
 to a subsequence if necessary):
\[
 A_p(\epsilon_i) = \displaystyle {\lim_{n\rightarrow \infty}} 
A_c(\varphi_n, D_{\epsilon_i}(p)),\qquad  K_p(\epsilon_i) = 
\displaystyle {\lim_{n\rightarrow \infty}} 
K(\varphi_n, D_{\epsilon_i}(p)), \qquad \forall i \in {\bf N},
\] 
and 
\[ A_p =   \displaystyle {\lim_{n\rightarrow \infty}} 
A_c(\epsilon_i),\qquad  K_p = \displaystyle {\lim_{n\rightarrow \infty}} K(\epsilon_i).
\]

For any  fixed $i$,  then $ \delta_n \cdot r < \epsilon_i $ if
$n$ is large enough. Equation~(\ref{eq:blow:bubb2}) then implies:
\[
A_c(\phi_n, D_{r}(p) ) = A_c(\varphi_n, D_{\delta_n \cdot r}(p))
 \leq A_c(\varphi_n,D_{\epsilon_i}(p)). 
\]

Taking the limit on both sides as $ n\rightarrow \infty,$ the result is:
\[
\displaystyle{\lim_{n \rightarrow \infty}} A_c(\phi_n, D_{r}) \leq
 A_p(\epsilon_i), \qquad \forall i \in {\bf N}.
\]

Taking limit on both sides as $ i \rightarrow \infty$, 
\[
\displaystyle{\lim_{n \rightarrow \infty}} A_c(\phi_n, D_{r}) \leq A_p.\]

Similarly, 
\[
\displaystyle{\lim_{n \rightarrow \infty}} K_c(\phi_n, D_{r}) \leq K_p.
\]

This implies that $ m \leq \sqrt{\frac{K_p \cdot A_p} {4 \pi^2}}.\; $
Applying Theorem~\ref{theo:lo:weak} for $\{\phi_n\}$ in $D_r$, we have: 
\[
\begin{array}{ccc} \displaystyle{\lim_{n \rightarrow \infty}} A_c(\phi_n, D_{r})& = &
   A_c(\phi_0, D_{r} \setminus \{q_1,q_2,\cdots, q_m\}) + \sum_{i=1}^{m} A_{q_i} \\
\displaystyle{\lim_{n \rightarrow \infty}} 
K(\phi_n, D_{r}) & \geq & K_c(\phi_0, 
D_{r} \setminus \{q_1,q_2,\cdots, q_m\}) + \sum_{i=1}^{m} K_{q_i}.\end{array}
\]

Thus,
\[
\begin{array}{ccc} A_p & \geq &  A_c(\phi_0, D_{r} \setminus \{q_1,q_2,\cdots, q_m\})
 + \sum_{i=1}^{m} A_{q_i} \\ K_p & \geq &  
K(\phi_0, D_{r} \setminus \{q_1,q_2,\cdots, q_m\}) + 
\sum_{i=1}^{m} K_{q_i}.\end{array}
\]

Let $r \rightarrow \infty,$ then:
\[
\begin{array}{ccc} A_p & \geq &  A_c(\phi_0, S^2 \setminus \{\infty, q_1,q_2,\cdots, q_m\})
 + \sum_{i=1}^{m} A_{q_i} \\ K_p & \geq &  
K(\phi_0, S^2\setminus \{\infty,q_1,q_2,\cdots, q_m\}) + \sum_{i=1}^{m} K_{q_i}.
\end{array}
\]

If a vanishing case occurs in  $D_{2^j}(0)$ for some $j,$
 then it occurs in any  disk $D_{2^i}(0).\;$ Observe the following
inequalities:
\[
 \begin{array}{ccl}  2 \pi \cdot e^{\phi_{n}(0) } & \geq & \int_0^{2 \pi}
 e^{\phi_{n}(\cos \theta,\sin \theta) } d\,\theta  \\
   & =  &   \int_0^{2 \pi} e^{\varphi_{n}(\delta_{n}\cos \theta,\delta_{n}\sin \theta) }
 \delta_{n} \cdot d\,\theta = \varepsilon.\end{array}
\]

The first inequality holds true because of equation~(\ref{eq:blow0}). The last two
equalities holds true because of  equation~(\ref{eq:blow4}) and~(\ref{eq:blow5}).  
According to lemma~\ref{lem:regular} ( in a vanishing case),
 the following two statements hold true: (1) $p$ is a bubble point of
 $\{\phi_{n}\}$;  (2) there exists at least one bubble point in the unit circle. Thus,
 in a vanishing case,  $m \geq 2.\;$  \\

QED.

\section{Geometrical Consequence}

\subsection{Theorem of weak convergence} 
A Riemannian metric is said to be a ``limit metric'' if
it is a weak limit of a sequence of Riemannian metrics
 in $H^{2,2}(\Omega).\;$ Lemma 6 implies that a limit metric
vanishes at one point if and only if it vanishes everywhere in its domain.
For the convenience of notations, we add the $ ``0'' $ metric into 
$H^{2,2}(\Omega)$ and the  resulting space is denoted by ${\hat{H}}^{2,2}(\Omega).\; $
Assume the following:
\[ K(0, \Omega) = A(0, \Omega) = 0, \qquad \mbox{ for any sub-domain}\, \Omega.\]
A sequence of Riemannian metrics $\{g_n\} \in H^{2,2}(\Omega)$
 weakly converges to a limit
metric $g_0$ in ${\hat{H}}^{2,2}_{loc}(\Omega) $ if and only if one of the following two 
alternatives holds true (mutually exclusive):
\begin{enumerate}
\item (Vanishing case). If $g_0 \equiv 0,$ then $g_n \rightarrow 0$ 
everywhere.
\item (Non-vanishing case). If $g_0 \neq 0, $ then
  $ \varphi_n \rightharpoonup \varphi_0 \; {\rm in}\; H^{2,2}_{loc}(\Omega),$
where $ g_n = e^{ 2\varphi_n} g_{bk}, g_0 = e^{ 2 \varphi_0} g_{bk}; $
and $g_{bk}$ is a smooth background metric in $\Omega.\;$
\end{enumerate}
We are now ready to re-state the theorem~\ref{theo:lo:weak} in
a geometric  context:\\

\noindent {\bf Theorem~\ref{theo:lo:weak}$'$.}
{\it  Let $\{g_n,\; n\in {\bf N}\}$ be a sequence of metrics
 with a finite area $C_1$ and energy $C_2 $ in a coordinate disk $D.\;$
 There exists a  subsequence $\{g_{n_j},j \in {\bf N}\}$ of $\{g_n\},$ 
 a finite number of bubble points
 $\{p_1,p_2,\cdots,p_m\}( 0\leq m \leq \sqrt{\frac{C_1\cdot C_2}{4\pi^2}} )$
 with respect to $\{g_{n_j}, j \in {\bf N}\},$  and a limit metric
 $g_0$ in ${\hat{H}}^{2,2}(D \setminus \{p_1,p_2,\cdots,p_m) $
such that:  
\[
 g_{n_j} \rightharpoonup g_0\;{\rm  in}\; {\hat{H}}^{2,2}_{loc}
(D \setminus \{p_1,p_2,\cdots,p_m\}).\]

If the amount of area and energy concentrations of $\{g_{n_j}\}$
 at each point $p_i$ are $A_{p_i}$ and $K_{p_i},$  then:
\begin{eqnarray}
 \displaystyle{\lim_{j\rightarrow \infty}} A(g_{n_j},D) 
&  = & A(g_0,D\setminus\{p_1,p_2,\cdots,p_m\} ) + \sum_{i=1}^{m} A_{p_i}  \\
  \displaystyle{\lim_{j\rightarrow \infty}} K(g_{n_j},D) & 
\geq & K(g_0,D\setminus\{p_1,p_2,\cdots,p_m\} ) +  
\sum_{i=1}^{m} K_{p_i}. \end{eqnarray}}\\

\noindent {\bf Proof.}  Re-write
 the sequence of metrics in a fixed coordinate system as:
\[  g_n = e^{2 \varphi_n} (d\,x^2 + d\,y^2). \]
Thus, $\{\varphi_n, n\in {\bf N}\}$ is a sequence of metrics with finite area
$C_1$ and energy $C_2.\;$ The rest of the proof is a direct translation
of the proof of  theorem~\ref{theo:lo:weak}
on p.~\pageref{theo:lo:weak}. QED.

\begin{theo} Let $\{g_n,\; n\in {\bf N}\}$ be a sequence of Riemannian metrics
in $M$ ($ M$ is any open surface)
 with a finite  area $C_1$ and energy $C_2.\;$ There exists a 
subsequence of $\{g_n \}, $ a finite number of 
bubble  points 
$\{p_1,p_2,\cdots,p_m\} ( 0\leq m \leq \sqrt{\frac{C_1\cdot C_2} {4 \pi^2}})$
with respect to this subsequence, and a limit metric $g_0$
 such that:
\[ g_n \rightharpoonup g_0 \;{\rm in}\;
 {\hat{H}}^{2,2}_{loc}(M \setminus \{p_1,p_2,\cdots,p_m\}).
   \]
If the amount of area and energy concentrations at each point $p_i$
are $A_{p_i}$ and $K_{p_i}$ , then:

\begin{eqnarray}
 \displaystyle{\lim_{j\rightarrow \infty}} A(g_{n_j},M) & 
 = & A(g_0,M\setminus\{p_1,p_2,\cdots,p_m\} ) + \sum_{i=1}^{m} A_{p_i} \\
  \displaystyle{\lim_{j\rightarrow \infty}} K(g_{n_j},M) 
& \geq & K(g_0,M\setminus\{p_1,p_2,\cdots,p_m\} ) +  \sum_{i=1}^{m} K_{p_i}.
\end{eqnarray}
\end{theo}

\noindent {\bf Proof}:
Let $\{U_1, U_2, \cdots, U_n,\cdots\}$ be a locally finite covering
 of $M$ where each $U_j$ is a  coordinate disk.
Consider the restrictions of the sequence of metrics $\{ g_n\}$ 
in each $U_j.\;$ These metrics  have a finite area $C_1$
and  a finite energy  $C_2.\;$ Apply theorem~\ref{theo:lo:weak} successively to
metrics in each coordinate disk. In $U_1,$  there exists a subsequence $\{g_{1n},
n\in {\bf N}\} $ of $\{g_n\},$  a finite set of bubble points 
$S_1 = \{q_{11},q_{12},\cdots,q_{1m_1}\}
(0\leq m_1 \leq\sqrt{\frac{C_1\cdot C_2}{4\pi^2}} )$
 with respect to this subsequence, and a limit metric $h_1 $ in
${\hat{H}}^{2,2}_{loc}(U_1\setminus S_1)$ such that:
\[ g_{1n} \rightharpoonup h_1 \;{\rm in}\; {\hat{H}}^{2,2}_{loc} (U_1 \setminus S_1).\]

Consider this subsequence $\{g_{1n}\}$ in $U_2.\;$ There exists a subsequence 
$\{g_{2n}, n\in {\bf N}\}$ of $ \{g_{1n}, n\in {\bf N}\},$ 
 a finite set of bubble points 
$S_2 = \{q_{21},q_{22},\cdots,q_{2m_2}\}
(0 \leq m_2 \leq\sqrt{\frac{C_1\cdot C_2}{4\pi^2}} )$ with respect to this subsequence, 
and a limit metric $h_2$ in ${\hat{H}}^{2,2}_{loc}(U_1 \setminus S_1) $ such that:
\[
 g_{2n} \rightharpoonup h_2 \;{\rm in}\;  {\hat{H}}^{2,2}_{loc} (U_2 \setminus S_2).\]

In general, if $\{g_{jn}\}$ had been defined in each coordinate disk $U_i ( i \leq j),$ 
 we can select a subsequence $\{g_{(j+1)n}\}$ of $\{g_{jn}\}$
in $U_{j+1}$ 
 so that there is a finite number of bubble points 
$S_{j+1}= \{q_{(j+1)1},q_{(j+1)2},\cdots,q_{(j+1)m_{j+1}}\}
(0 \leq m_{j+1} \leq \sqrt{\frac{C_1\cdot c_2}{4\pi^2}} )$
in $U_{j+1}$ with respect to this subsequence, and a limit metric
$h_{j+1}$ in ${\hat{H}}^{2,2}_{loc}(U_{j+1}\setminus S_{j+1})$
 such that:
\[
 g_{(j+1)n} \rightharpoonup h_{j+1} \;{\rm in}\;
  {\hat{H}}^{2,2}_{loc} (U_{j+1} \setminus S_{j+1}).
\]
Consider the diagonal  subsequence $\{g_{nn}, n\in {\bf N}\}.\;$ 
In each coordinate disk $U_j (\forall j)$,
the following holds true:
\[ g_{nn} \rightharpoonup h_{j} \; {\rm in} \;
{\hat{H}}^{2,2}_{loc} (U_{j} \setminus S_{j}).\]

This set of limit metrics $\{h_j\}$ then defines a limit metric
 $g_0$ in $H^{2,2}_{loc} (M \setminus (\bigcup_{j} S_j))$
by:
\[ g_0(p)  = h_j (p),\qquad \forall\; p\, \in U_j.\]

This metrics is well defined since  $h_i \equiv h_j$ on $U_i \bigcap U_j$ if 
$U_i \bigcap U_j \neq \emptyset.\;$
Thus,
\[
g_{nn} \rightharpoonup g_0 \;{\rm in}\; H^{2,2}_{loc} (M \setminus (\bigcup_{j} S_j)).
\] 

The cardinality of the set $\bigcup_{j} S_j $ must be
 bounded by $ \sqrt{\frac{C_1 \cdot C_2}{4 \pi^2}} $
according to the proof of proposition 2. \\
 
Re-Label this subsequence as $\{g_n\}.\;$
For any two pair of coordinate
disks $U_i,\; U_j$ where $ U_i \bigcap U_j \not = \emptyset,$  we have:
\begin{eqnarray}
 \displaystyle{\lim_{n\rightarrow \infty}} A(g_n,U_j \bigcup U_k) 
&  = & A(g_0, (U_j \bigcup U_k)\setminus(S_j \bigcup S_k) )
 + \sum_{ p \in S_j \bigcup S_k} A_{p},
 \\
  \displaystyle{\lim_{n\rightarrow \infty}} K(g_n,U_j \bigcup U_k) 
& \geq & K(g_0,(U_j \bigcup U_k)\setminus (S_j \bigcup S_k))
 +  \sum_{ p \in S_j \bigcup S_k} K_{p}. 
\end{eqnarray}
These two formulas  can be readily generalized to  any number of coordinate disks:
\begin{eqnarray}
 \displaystyle{\lim_{n\rightarrow \infty}} A(g_n, \bigcup_{k} U_k)
 &  = & A(g_0, \bigcup_{k} U_k \setminus ( \bigcup_k S_k) )
 + \sum_{ p \in  \bigcup_k S_k} A_{p},
 \\
  \displaystyle{\lim_{n\rightarrow \infty}} K(g_n, \bigcup_k U_k) 
& \geq & K(g_0, \bigcup_k U_k \setminus ( \bigcup_k S_k)) 
+  \sum_{ p \in \bigcup_k S_k} K_{p}. 
\end{eqnarray}
Observed that $M= \bigcup_{k} U_s $ and $ \bigcup_k S_k = \{p_1,p_2,\cdots,p_m\}.\;$
QED.\\

\subsection {Blowing up procedure and tenuously connected sum}
Let us re-state theorem~\ref{theo:lo:bubble} in the geometric context.
 \\

\noindent {\bf Theorem ~\ref{theo:lo:bubble}$'$} (Bubbles on bubbles). {\it
 Let $\{g_n, n\in{\bf N}\}$ be a sequence of metrics
in $D$  with a finite area $C_1$ and energy $C_2.\;$ Suppose that $ p=0$ is
the only bubble point in $D$ with area concentration $A_p$ and energy
concentration $K_p.\;$  Fix a local $z-$coordinate system centered at $p$
and a scaling constant $\varepsilon.\;$ If $\varepsilon$ is small enough,
we can re-normalize the sequence of metrics by $\tilde{g}_n(z)
 = g_n(\epsilon_n \cdot z + z(p_n))$
where $\{\epsilon_n \searrow 0 \}$ is uniquely
determined by the scaling constant $\varepsilon;$ 
where  $p_n \rightarrow p$  is the supremum of mass $g_n$ in  $D.\;$
There then exists a subsequence of $\{g_n\}$,  a finite number of bubble points 
$\{q_1,q_2,\cdots,q_m\} (0\leq m \leq \sqrt{\frac{A_p \cdot K_p} {4\pi^2}})$
in $S^2 \setminus \{\infty\}$ with respect to the corresponding
subsequence of $\{\tilde{g}_n\},$ and a limit metric $\tilde{g}_0$
 in ${\hat{H}}^{2,2}(S^{2}\setminus \{\infty, q_1,q_2, \cdots,q_m\})$
such that:
\[\tilde{g}_n \rightharpoonup \tilde{g}_0\; {\rm in }\; {\hat{H}}^{2,2}_{loc} 
(S^{2} \setminus \{\infty, q_1,q_2, \cdots,q_m\} ).\]

If the  amount of  area and energy concentrations of $\{\tilde{g}_n\}$
 at each $q_i $   are $A_{q_i}$ and $K_{q_i}$  respectively, 
then:
\begin{eqnarray}  
 A_p & \geq &  A(\tilde{g}_0,S^2 \setminus\{\infty,q_1,q_2,\cdots, q_m\}) + \sum_{i=1}^{m} A_{q_i}, \\
K_p & \geq &  K(\tilde{g}_0,S^2 \setminus\{\infty,q_1,q_2,\cdots, q_m\}) + \sum_{i=1}^{m} K_{q_i}.
\end{eqnarray}
}\\

Let us review the steps taken in the proof of theorem~\ref{theo:lo:bubble}.
For convenience, we use a complex notation. The metric can be expressed as:
\[
g_n(z) = e^{2 \varphi_n(z)} |d\,z|^2, \qquad \forall n \in {\bf N}.
\]

The first step is to move the supremum of mass of the metric $g_n$ to the center 
of the coordinate system. If $\{p_n\}$ is such a sequence of points,
then  define
\[
  \tilde{g}_n (z) = g_n(z + z(p_n)),\qquad z \in  D_{1}(p).
\] 

Re-Label $\{\tilde{g}_n \}$ as $\{g_n\}.\;$ The supremum of
the metric $g_n$ is now at $p,\;\forall n \in {\bf N}.\;$ \\

Choose a small positive number $\varepsilon < \epsilon_0 $(as in 
corollary~\ref{cor:blow}) as a filter.
Following the proof of theorem 3,  there then exists
 a number  $ r_{1} > 0 $ such that if  $ n $ is large enough, we have: 
\begin{equation} \displaystyle{\max_{0\leq\theta\leq 2\pi}}
 ( {\varphi_{n}}(r_{1} \cos \theta,r_1 \sin \theta) + \ln r_{1}) 
 \ll \varepsilon, \qquad \forall \;n > N, 
\label{eq: glo:blow1}
\end{equation}
or the length of this circle at $ |z| = r_1$ is very small:
\begin{equation}
 \int_0^{2 \pi} e^{\varphi_n(r_1 \cdot z)}\, r_1\, d\,\theta 
\ll \varepsilon, \qquad n > {\bf N}.
\label{eq:glo:blow2}
\end{equation}

Following corollary~\ref{cor:blow},there exists  $\delta_{n} > 0 $ such that: 
\begin{equation} \int_0^{2 \pi} e^{\varphi_n(r \cdot z)}\, r\, d\,\theta 
 < \varepsilon,\qquad \forall\;  r_{1} \geq r \geq {\delta_{n}}, 
\label{eq:glo:blow3} \end{equation}
and 
\begin{equation}
\int_0^{2 \pi} e^{\varphi_n (\delta_n \cdot z)}\, \delta_n \,d\,\theta  = \varepsilon.  
\label{eq:glo:blow4}
\end{equation}

The circle $ |z| = \delta_n$ is the first circle for which the
metric $g_n$  has a length of $\varepsilon\;$ beyond
a  thin neck.  Hence, the set of concentric circles $\{|z| = \delta_n \}$
is uniquely determined by the filter size $\varepsilon $ once the local coordinate
system is picked.  Define a sequence of conformal parameter functions as:
\[ \phi_n(z) = \varphi_n( \delta_n \cdot z) + \ln \delta_n, \qquad n \in {\bf N}.\]

Thus re-normalize the original sequence of  metrics as
\[
\tilde{g}_n(z) = e^{2 \phi_n(z)} |d\,z|^2 =  g_n( \delta_n \cdot z), \qquad n \in {\bf N}.
\]  

Theorem~\ref{theo:lo:bubble} then asserts that we could 
choose a subsequence $\{\varphi_{n_i}, i\in {\bf N}\}$ of
$\{\varphi_n, n\in {\bf N}\},$  a finite number of bubble
points $\{q_1,q_2,\cdots,q_m\} ( 0\leq m \leq \sqrt{\frac{A_p\cdot K_p}{4\pi^2}})$
such that one of the following two alternatives holds true:
\begin{enumerate}
\item $\phi_{n_j}( z) 
   \rightarrow -\infty $ in $ S^2\setminus \{\infty,q_1,q_2,\cdots,q_m\}.$ 
\item There exists a metric $\phi_0 \in H^{2,2}_{loc}( S^2\setminus
\{\infty,q_1,q_2,\cdots,q_m\})$  such that
\[
\phi_{n_j}( z)  
\rightharpoonup \phi_0, \;{\rm in}\;  H^{2,2}_{loc}( S^2\setminus
\{\infty,q_1,q_2,\cdots,q_m\}). 
 \]
\end{enumerate}

Define $\tilde{g}_0 \equiv 0 $ in a vanishing case; and define
$\tilde{g}_0 = e^{ 2 \phi_0} |d\,z|^2 $ in a non-vanishing case. Thus,
\[ \tilde{g}_{n_j} \rightharpoonup \tilde{g}_0\;{\rm in}\;   {\hat{H}}^{2,2}_{loc}( S^2\setminus
\{\infty,q_1,q_2,\cdots,q_m\}).  \]

Moreover,
\begin{eqnarray}  
 \tau_p   & = & A_p - A(g_0,S^2 \setminus\{\infty,q_1,q_2,\cdots, q_m\})
 - \sum_{i=1}^{m} A_{q_i} \geq 0,\\
K_p & \geq &  K(g_0,S^2 \setminus\{\infty,q_1,q_2,\cdots, q_m\})
 + \sum_{i=1}^{m} K_{q_i},
\end{eqnarray}
 
where $\tau_p$ denote the amount of area lost in the neck during the 
re-normalization (blowing up) process.\\

Choose $r_2 $ big enough, so that $\{q_1,q_2,\cdots,q_m\} 
\subset D_{r_2}.\;$ Consider the cylinder bounded
by the two concentric circles $|z| = r_1$ and $|z| = r_2 \cdot \delta_n.\;$
This cylinder is called the ``neck'' of the blowing up process. 
By definition, the length of a circle in this cylinder is bounded
above by $\varepsilon.\; $ As $n\rightarrow
\infty,$  the conformal distance between the two boundary circles approaches $\infty,$ while
part of the interior of the neck collapses into a line. 
The collapsing can occur  either by keeping the
 scalar curvature point-wisely bounded, or by keeping  the diameter of the neck 
bounded. Denoted this neck by $Neck(r_1,r_2).\;$ We can shrink
the size of the neck by letting $r_1 \rightarrow 0$ and $r_2 \rightarrow \infty.$\\ 

This blowing up procedure or the re-normalization procedure depends only on the
filter size $\varepsilon> 0 $ once a coordinate system is fixed. 
Suppose that $g_0$ is a limit  metric in ${\hat{H}}^{2,2}(D\setminus\{p\})$ such that:

\[ g_n \rightharpoonup g_0\; {\rm in}\;{\hat{H}}^{2,2}_{loc}(D\setminus\{p\}). \]
The surface $(g_0, D\setminus{p})$ and 
$(\tilde{g}_0, S^{2} \setminus \{\infty,q_1,q_2,\cdots,q_m\}) $ 
are called tenuously connected at $p$ and at $z=\infty.\;$
If $\tau_p = 0,$   the connected sum is then efficient. Otherwise 
the connected sum is inefficient.\\

The following is an example of a ``tenuously connected sum'' of two
Riemannian metrics or surfaces:\\

\noindent {\bf Example 2.} {\it We first construct 
a metric in a disk, where the boundary curve is a closed geodesic. We can make
the length of the boundary approaches $0,$ while keeping the area and energy
finite (see Figure~\ref{fg: closed geod.} below).
 The following is a sketch of the construction. Suppose that
$g = e^{2\varphi} |d\,z^2|$ is a rotationally symmetric metric defined in
 ${\bf R^2}$ (real plane) such that: \label{example}
\[
  \varphi(r) = -\ln r - \beta \cdot \ln (\ln r),\qquad \forall\,r > 2, 
\]
where $ {1\over 2} < \beta < {3\over 2}.\;$
Let $\epsilon_n = {{\beta} \over{\ln n}},\delta_n = {{\beta} \over{\ln^2 n}},$ and
 $ T_n = n + \ln n.\;$ Define a sequence of metrics $\{\varphi_n\} $ in $D_{T_n}$
as the following:
\[
\varphi_n(r) =\left\{\begin{array}{ll} \varphi(r), &\mbox{when $r\leq n$}\\
                                      \varphi(n)+\ln n -\ln r -\epsilon_n (r-n) +{1\over 2}\delta_n
    (r-n)^2, &\mbox{when $n \leq r \leq T_n$} \end{array} \right.
\]

It is straightforward to prove that $|z| = T_n$ is a closed geodesic of $\varphi_n$
and 
\[  \displaystyle{\lim_{n\rightarrow \infty}} \varphi_n(r) = \varphi(r), \qquad \mbox{ if $r$ is finite,}\]
and
\[
 \displaystyle{\lim_{n\rightarrow \infty}} E_c(\varphi_n,D_{T_n}) = E_c(\varphi,{\bf R^2}),\qquad
 \displaystyle{\lim_{n\rightarrow \infty}} A_c(\varphi_n,D_{T_n}) = A_c(\varphi,{\bf R^2}).
\]

Gluing two identical copy of $\varphi_n$ along
the  curve $|z| = T_n,$ we obtain a metric $g_n$ in $S^2.\;$
Clearly, $\{g_n\}$ has finite energy and area, and it weakly 
  converges to $g$ everywhere except at
 near $z=\infty.\;$ If we blow up the sequence near $z=\infty$, we obtain a 
 new sequence  of metrics
which locally weakly converges to $g$ except at $z=0.\;$ 
Re-label this metric as $g_{\infty}.\;$  Then the limit tree structure of the 
weak limit of $(S^2, g_n)$  consists of a root vertex $g$ and a child vertex
$g_{\infty}.\;$ The corresponding metrics $g$ and $ g_{\infty}$ at 
the two nodes are tenuously connected. This sequence of metrics, clearly has
no convergent subsequence in the elementary sense, even up to the
M\"obious group. Using a similar mechanism,  we could construct examples
of a sequence of metrics which demonstrates a more sophisticated pattern 
of limit tree structures.
}\\

\begin{figure}
\centerline{\psfig{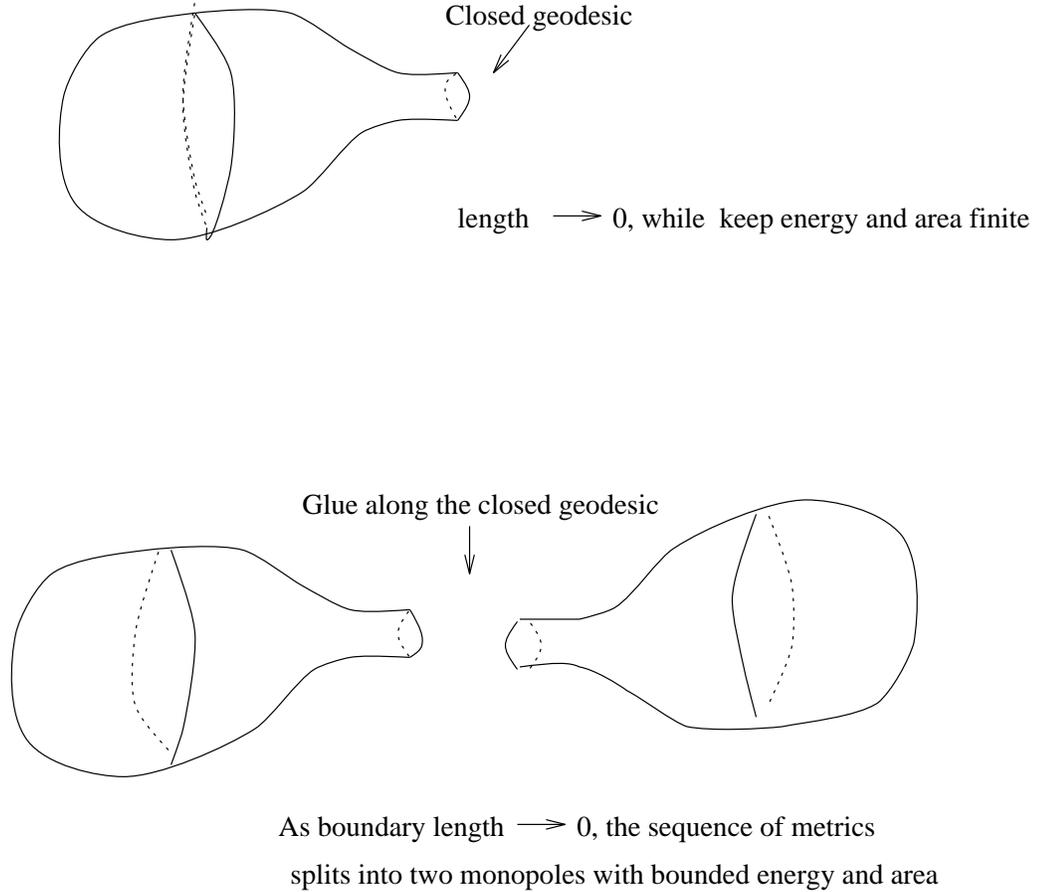}}
\caption{Surface bounded with a small closed geodesic}
\label{fg: closed geod.}
\end{figure}

\begin{lem} Let  $g $ be a metric in  $D\setminus\{p\}$ with a finite energy $C_2.\;$
Suppose $|\partial D_1|_g = \epsilon $ where $D_1 \subset D\setminus\{p\}.\;$
There exists a constant $C_{\epsilon} > 0 $ (independent of
metric $g$) such that:
\[ A(g, D\setminus \{p\}) = \int_{D\setminus \{p\}} d\,g > C_{\epsilon}.\]
\label{lem:leastarea}
\end{lem}
{\bf Proof.} If the lemma is false, there then exists a
sequence of metrics $\{g_n\}$ in $D\setminus\{p\}$
such that:
\[K(g_n,D\setminus\{p\}) < C_2,\qquad A(g_n,D\setminus\{p\}) 
\rightarrow 0,\qquad |\partial D_1|_{g_n} = \epsilon>0.\]

According to theorem 5, there exists a subsequence
of $\{g_{n}\},$  a finite bubble points $\{p_1, p_2, \cdots,p_m\} (m \geq 0)$, and
a limit metric $g_0$ in ${\hat{H}}^{2,2}(D\setminus\{p,p_1,p_2,\cdots,p_m\})$
 such that:
\[ g_n \rightharpoonup g_0\;{\rm in}\; {\hat{H}}^{2,2}_{loc}(D\setminus\{p,p_1,p_2,\cdots,p_m\}).
\]

Thus, $m = 0 $ since the product of area and energy of this subsequence approaches
$0.\;$ Moreover, $g_0 \equiv 0 $ in $ D\setminus\{p\} $ since total
area approaches $0.\;$
This is also impossible since $ |\partial D_1|_{g_n} = \epsilon > 0.\;$
The lemma is then proved. QED.\\

\begin{prop} (continued from theorem 3$'$). $\tilde{g}_0$ is as defined in theorem 
3$'$. If $\tilde{g}_0\neq 0$ and $m=1,$ then 
$ \int\limits_{S^2 \setminus\{\infty,p\}} \,d\,\tilde{g}_0
 > C_{\varepsilon} > 0,$ where $C_{\varepsilon}$ depends only 
on $C_1, C_2$ and the scaling constant $\varepsilon.$
\end{prop}
Proof: If $\tilde{g}_0 \neq 0 $ and $m=1, $ then $p$ must be the only bubble point
and $|\partial D_1 |_{g_0} = \epsilon > 0.\;$ The proposition
then follows the previous lemma. QED.\\

\section{Bubble tree}

\noindent {\bf Theorem A.}\label{th:theorem A} {\it The limit of any locally weakly convergent 
sequence of metrics $\;\{g_k,\;k \in \bf {N}\}\;$  $\in\;$  
${\cal{S}}(g_{0},C_1,C_2,\Omega),\;$
  consists of the following 4 objects: (1) A 
finite, rooted tree $T$, possibly reduced to just the base vertex $f.$  
(2) The base vertex $f \in T $ is a
limit metric in $\Omega$ with a finite number of bubble points $\{p_i\}$
deleted; the edges emanating from the base vertex is $\{p_i\};$
there are three masses associated with each edge: the area
concentration $a_i$, energy concentration $e_i$ and
area loss during the blowing up process $\tau_i$ ($ a_i \cdot e_i \geq 4\pi^2 $). (3)
Any other vertex $f_s$ is a  limit metric defined
on $S^{2}\setminus\{\infty,p_{si}\};$ the edges emanating
from this vertexes are $\{p_{si}\}$; and there are three masses associated
with each edge:  the area concentration $a_{si}$, energy concentration $e_{si}$ and
area lost during the blowing up process $\tau_{si}.\;$
(4) For each pair of vertexes $f_{s_1} $ and $ f_{s_2}$ bounding 
an common edge in $T$, they are tenuously connected at the pair of respective 
singular points. The connected sum is efficient if the area loss associated with
that edge is $0.\;$
If the tree $T$ consists of
only the base vertex $f,\;$the sequence of metrics $\{g_k\}$ is then said
to have a weak convergent limit in the elementary sense (up to the M\"obious group).
The number of vertexes whose valence $\neq 2$, is bounded from
above  $ (\leq \sqrt{ C_1\cdot C_2}).\;$
 The depth of the tree is also finite in a reasonable sense. } \\

\noindent {\bf Proof of Theorem A}.
The tree structure is constructed from a sequence of metrics $\{g_n\}$
 in $\Omega$ as follows (see Figure~\ref{fg:closed} on p.~\pageref{fg:closed}):
 First, choose a scaling constant $\varepsilon_0$ as
a filter for re-normalization process.  The $\{g_n\}$ locally weakly converges to 
$f_0$ on $\Omega \setminus\{p_1,p_2,\cdots,p_m\} $ except a finite number of
bubble points $\{p_1,p_2,\cdots,p_m \}.\;$ The base vertex of the
tree is the metric $f_0$, which we re-labeled as $f,$
and the edges emanating from the base vertex  are the points $\{p_i\}.\;$
Each edge has an energy mass $e_i$ and area mass $a_i,$ which are the energy and area
concentrations at the bubble point $p_i.\;$  For each $p_i$, the re-normalization
process gives a new sequence of metrics  $\{\tilde{g}_n\}$ 
which locally weakly converge to a metric $f_i$
in $S^2 \setminus\{\infty,p_{ij},j=1,2,\cdots,m_i\},$ with the
amount of energy and area concentrated at
each point $p_{ij}$ are $e_{ij}$ and $a_{ij}.\;$ We label each edge as $(p_i,e_i,a_i, \tau_i)$
where $\tau_i$ represent the amount of area lost at the bubble point $p_i$
during the blowing up process.   If $\tau_i = 0, $  the blowing
up process is then efficient.  The edge $(p_i,e_i, a_i, \tau_i)$ ($e_i \cdot a_i \geq 4\pi^2 $ according to lemma 2) terminates 
at the vertex $f_i$ which, in turn, is the source of new edges $\{p_{ij}\}, $
 and so on. \\

At each vertex $f_I = f_{i_1\cdots i_{k-1}i_k}$ of the tree, use 
$S_I $ to denote all of the bubble points of this limit metric other than
the point $ z=\infty.\; $ If $f_I$ is not the base vertex, it must have a
 parent vertex   $f_{I'} = f_{i_1 \cdots i_{k-1}}$. The surface
$(f_{I'}, S^2\setminus\{\infty, S_{I'}\}) $  
( or $ (f, M\setminus\{p_1,p_2,\cdots,p_m\} $ if $I' = \emptyset $)
 is tenuously connected to
 $( f_{I},S^2\setminus\{\infty, S_{I'}\}).\; $ If there is area loss
during the blowing up process ($\tau_{I'} \not = 0$), 
the connected sum is inefficient.
\\
 
 Each vertex $f_I$ has a special property: if it vanishes in any point
in its domain, then it vanishes everywhere in its domain. In the
case when $f_I \equiv 0,$ we call this a ghost vertex. At each 
ghost vertex other than the base vertex, there exists at least two edges
emanating from it. In other other words, the metric has at least two bubble
points.\\

The ghost vertex does appear, as seen in example 3 below. 
However, there exists at most a finite number
of ghost vertexes. Otherwise, consider all the vertexes in the tree that
have at least two edges emanating from them. These vertexes must
be infinitely many  since every ghost vertex has at least two edges emanating 
from it. There exists an infinite number of edges
where no two edges  belong to the same branch of the tree.
 Re-Labeling these edges if necessary, we may assume that these 
edges are $\{(q_i,e_i,a_i, \tau_i), i \in {\bf N}\}$
 where $a_i \cdot e_i \geq 4 \pi^2.\;$ Therefore,
\[ C_1 \geq \sum_i a_i, \qquad C_2\geq \sum_i e_i. \]
Thus,
\[ C_1 \cdot C_2 \geq \sum_i a_i \cdot e_i = \sum_i 4\pi^2.\]
The last inequality implies that the number of these vertexes (include
all the ghost vertexes)  must be finite. \\

For any other vertex which has only one edge emanating from
it, proposition 3 implies that  the area of such a vertex is bounded
below by  a positive constant $C_{\epsilon}$, 
which depends only on $C_1,C_2$ and  the scaling constant $\varepsilon.\;$
The  number of these vertexes is finite as well. \\

Therefore, the limit tree has only a finite depth. If we reduce the size
 of the filter, a new vertex might be  inserted into the tree 
structure. However, these
new vertexes have only one edge emanating from it. The underlying surface
is $S^2$ with two opposite points deleted. QED. \\

\noindent{\bf Example 3}. {\it Let $f= (z-1)(z-2) \cdots(z-m)$ be a 
holomorphic function. Choose a simply connected domain $\Omega$
which contains all zero points of $f$ but no zero points of $f'(z).\;$
Thus, $g_n(z) = { 4\cdot n^2 \cdot |f'|^2 \over{ (1 + n^2 \cdot |f|^2 )^2}} \cdot |d\,z|^2$
is a sequence of metrics well defined in $\Omega $ with finite area
and energy (bounded above by $4\pi \cdot m$).  Clearly, $g_n$
weakly converges to $0$ everywhere except at $z=1,2,\cdots,m.\;$ At each
point $z=k$, a renormalized sequence of metrics weakly converges to 
a metric in $S^2$ with curvature $1.\;$ Thus, the bubble tree
of $(g_n,\Omega)$ consists of 1 ghost base vertex and $m$ first
level vertexes, where each first level vertex represents a metric
with curvature 1 in $S^2$. }\\

\noindent {\bf Proof of Corollary B.} Suppose $ \{g_k,k \in {\bf N} \}$ is 
a sequence of metrics with finite area $C_1$ and energy $C_2$. 
If necessary, we pass to a subsequence so that the  weak limit of
this sequence  has a  bubble tree decomposition as
described in theorem A.  Consider a generic  pair of consecutive vertexes
$(f_I, f_{Ii})$ in the bubble tree, where $p_i$ is a bubble point of
$f_I$ and the re-normalized sequence of metrics at $p_{Ii}$ weakly converges
to $f_{Ii} $ except a few bubble points. Consider the ``neck'' of this blowing up process.
It is a cylinder where the length of each concentric circle 
is bounded above by the scaling constant  $\epsilon.\;$ We call this cylinder
a ``thin component.'' Now iterate thorough each pair of consecutive vertexes, and obtain
a finite number of ``thin'' components (See Figure~\ref{fg: thin-thick} below). The collection of ``thin
components'' is labeled  by $I_{thin}.\;$ For each fix $n$, remove
all of the ``thin components'' from $M.\;$ The resulting surface
is a disjoint union of  a finite number of connected components.
Each connected component is called a ``thick component.''  Label 
all of the ``thick'' components by $I_{thick}.\;$ Each thick component,
together with the restriction of $g_n$ on it, weakly converges to a surface
with a finite number of disks deleted. The thick component corresponding to the
base vertex is $\Omega$ with a few disk deleted. All of the rest of thick
components are $ S^2$ with a few disks deleted (When a thick component
corresponds to a ghost vertex in the tree decomposition, the limit
metric is $0$). The size of all deleted disks
could be shrinked to $0$ by shrinking the size of corresponding blowing
up ``neck.'' QED. \\

{\bf Acknowledgments.}  The author wishes to thank Professor J. Goodman
for carefully reading my manuscript and many helpful comments he made.
Thanks also go to Professor R. Schoen for many helpful
and stimulating conversations during the course of this work. 
The author also wishes to thank his advisor Professor E. Calabi
for his warm encouragement and continue supports during the past two years.\\

\begin{figure}
\centerline{\psfig{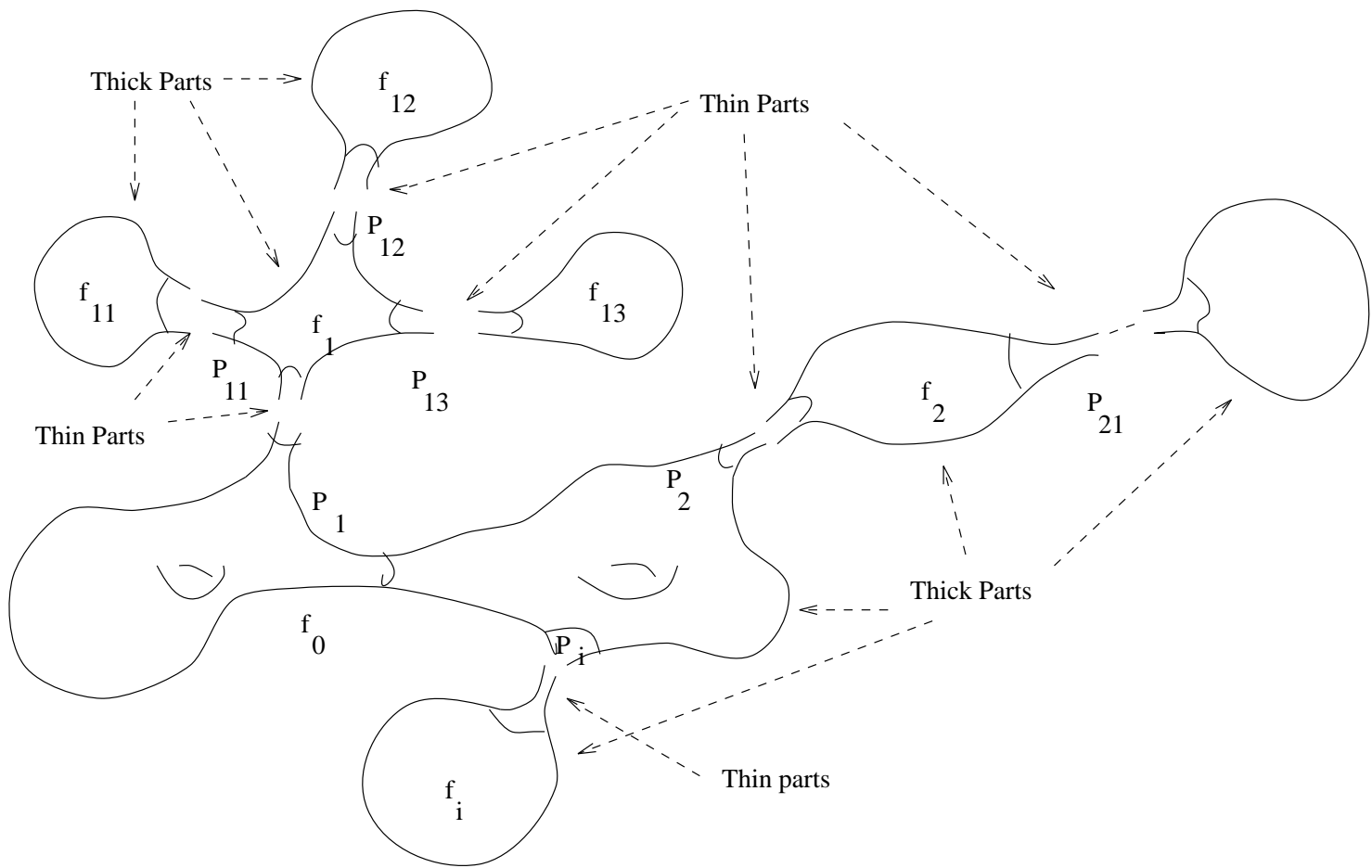}}
\caption{Thin-thick decompsition}
\label{fg: thin-thick}
\end{figure}

\nocite{Tshioya96}
\bibliography{test}

\end{document}